\documentclass[12pt]{amsart}
\usepackage[T1]{fontenc}
\usepackage[latin9]{inputenc}
\usepackage{geometry}
\usepackage{mathrsfs}
\usepackage{amstext}
\usepackage{amsthm}
\usepackage{amssymb}
\usepackage{bbm}
\usepackage[table]{xcolor}

\usepackage{algorithm}
\usepackage{algorithmic,eqparbox,array}

\usepackage{subcaption}

\usepackage[svgnames]{xcolor}
\usepackage[bookmarksnumbered=true]{hyperref} 
\hypersetup{
	colorlinks = true,
	linkcolor = Blue,
	anchorcolor = blue,
	citecolor = Green,
	filecolor = blue,
	urlcolor = FireBrick
}

\usepackage[normalem]{ulem}
\usepackage{aliascnt}
\usepackage[nameinlink]{cleveref}

\makeatletter
\numberwithin{equation}{section}
\numberwithin{figure}{section}

\theoremstyle{plain}
\newtheorem{theorem}{Theorem}[section]

\theoremstyle{definition}
\newaliascnt{definition}{theorem}

\aliascntresetthe{definition}
\crefname{definition}{Definition}{Definitions}

\theoremstyle{definition}
\newaliascnt{question}{theorem}

\aliascntresetthe{question}
\crefname{question}{Question}{Questions}

\theoremstyle{definition}
\newaliascnt{remark}{theorem}

\aliascntresetthe{remark}
\crefname{remark}{Remark}{Remarks}

\theoremstyle{plain}
\newaliascnt{corollary}{theorem}

\aliascntresetthe{corollary}
\crefname{corollary}{Corollary}{Corollaries}

\theoremstyle{plain}
\newaliascnt{lemma}{theorem}
\newtheorem{lemma}[lemma]{Lemma}
\aliascntresetthe{lemma}
\crefname{lemma}{Lemma}{Lemmas}

\theoremstyle{plain}
\newaliascnt{proposition}{theorem}

\aliascntresetthe{proposition}
\crefname{proposition}{Proposition}{Propositions}

\theoremstyle{definition}
\newaliascnt{example}{theorem}

\aliascntresetthe{example}
\crefname{example}{Example}{Examples}

\theoremstyle{definition}
\newaliascnt{assumption}{theorem}
\newtheorem{assumption}[assumption]{Assumption}
\aliascntresetthe{assumption}
\crefname{assumption}{Assumption}{Assumptions}

\theoremstyle{definition}
\newaliascnt{problem}{theorem}

\aliascntresetthe{problem}
\crefname{problem}{Problem}{Problems}

\newtheorem*{definition*}{Definition}
\newtheorem*{remark*}{Remark}
\newtheorem*{example*}{Example}

\renewenvironment{proof}[1][\proofname]{\medskip \noindent {\bfseries #1. }}{\hfill \qedsymbol\medskip}

\MakeRobust{\ref}
\newcommand{\labeltext}[2]{
\@bsphack
\csname phantomsection\endcsname
\def\@currentlabel{#1}{\label{#2}}
\@esphack
}

\usepackage{scalerel}
\usepackage[usestackEOL]{stackengine}
\def\dashint{\,\ThisStyle{\ensurestackMath{%
  \stackinset{c}{.2\LMpt}{c}{.5\LMpt}{\SavedStyle-}{\SavedStyle\phantom{\int}}}%
  \setbox0=\hbox{$\SavedStyle\int\,$}\kern-\wd0}\int}

\usepackage{enumitem}

\DeclareRobustCommand{\SkipTocEntry}[5]{}

\newcommand{\mR}{\mathbb{R}}   
\newcommand{\mC}{\mathbb{C}}   
\newcommand{\mN}{\mathbb{N}}   
\newcommand{\mS}{\mathbb{S}}   

\newcommand{\abs}[1]{\lvert #1 \rvert}  
\newcommand{\norm}[1]{\lVert #1 \rVert}  

\newcommand{\mI}{\mathcal{I}}
\newcommand{\mA}{\mathcal{A}}

\newcommand{\calS}{\mathcal{S}}
\newcommand{\calF}{\mathcal{F}}

\newcommand{\bfi}{\mathbf{i}}

\newcommand{\sfm}{\mathsf{m}}

\newcommand{\rmd}{\mathrm{d}}

\makeatother

\usepackage{babel}

\newcommand{\pl}{\partial}

\usepackage{cancel}

\usepackage{orcidlink}

\begin{document}

\title[An inverse problem for semilinear elliptic equations]{An inverse problem for semilinear elliptic equations with generalized Kerr-type nonlinearities}

\begin{sloppypar}

\begin{abstract}
We study the inverse problem of reconstructing the convex hull of unknown inclusions in semilinear elliptic equations with nonanalytic nonlinearities, by extending Ikehata's enclosure method to accommodate such nonlinear effects. 
To overcome the lack of the higher-order differentiability required for the Taylor-based construction, we construct an approximate solution based on the linearized equation through a first-iteration approximation.
Under suitable structural conditions on the nonlinearity, we establish a reconstruction result for the convex hull of the inclusion. 
The method applies to a broad class of generalized Kerr-type nonlinearities that need not be analytic, without requiring an explicit formula for the nonlinear term. We also address the numerical challenges arising from the exponential decay of the CGO test data and show how they limit the accuracy of an implementation in double-precision arithmetic.
\end{abstract}


\author[Kow]{Pu-Zhao Kow\,\orcidlink{0000-0002-2990-3591}$^{*}$}
\thanks{$^{*}$Department of Mathematical Sciences, National Chengchi University, Taipei 116, Taiwan (\texttt{pzkow@g.nccu.edu.tw}) \& Kanazawa University, Kakuma, 920-1192 Kanazawa, Japan. Kow was supported by the National Science and Technology Council of Taiwan, NSTC 112-2115-M-004-004-MY3 and NSTC 115-2628-M-004-001, and by the National Center for Theoretical Sciences of Taiwan.} 

\author[Kuan]{Rulin Kuan\,\orcidlink{0000-0003-0056-2380}$^{\dagger}$}
\thanks{$^{\dagger}$Department of Mathematics, National Cheng Kung University, Tainan 70101, Taiwan (\texttt{rkuan@ncku.edu.tw}). Kuan was supported by the National Science and Technology Council of Taiwan, NSTC 113-2115-M-006-005-MY2.}

\maketitle


\subsection*{Acknowledgments} 
\addtocontents{toc}{\SkipTocEntry} 

We would like to thank the anonymous readers for their comments and for pointing out several typos in the first draft on arXiv. 
We are also grateful to one of the anonymous referees for identifying a key issue in the numerical simulations that had escaped our attention during the previous rounds of revision. In response, we have added \Cref{sec:numerical} to present the numerical implementation and address this issue.

\subsection*{Code Availability} 
\addtocontents{toc}{\SkipTocEntry}

A MATLAB implementation of \Cref{alg:implementation} can be found in our Github repository\footnote{\url{https://github.com/puzhaokow1993/nonlinear_enclosure_method}}. 

\subsection*{AI declaration} 
\addtocontents{toc}{\SkipTocEntry}

Claude Opus 5.5 was used to assist in developing the theoretical analysis and numerical implementation presented in \Cref{sec:numerical}. In particular, it assisted in assessing the limitations of double-precision arithmetic in the design and MATLAB implementation of \Cref{alg:implementation}, as well as in addressing the associated theoretical challenges in establishing \Cref{lem:well-posed-relax} and \Cref{thm:approximation}. The authors assume full responsibility for the content of the manuscript.

\section{Introduction\label{sec:introduction}}

In this paper, we consider the inverse problem of reconstructing the convex hull of an unknown inclusion in a domain from boundary measurements, governed by the following class of semilinear elliptic equations. Let $\Omega \subset \mathbb{R}^n$ ($n \ge 2$) be a bounded $C^{\infty}$-domain. We study the problem:
\begin{subequations} \label{eq:main}
\begin{equation}
\Delta u_{f} + q(x,\abs{u_{f}})u_{f} = 0 \text{ in $\Omega$} ,\quad \left.u_{f}\right|_{\partial\Omega}=f. \label{eq:generalized}
\end{equation}
Here, the nonlinear term $q(x,\abs{u})$ in our model is non-analytic in $u$, and is designed to represent a range of physically relevant responses. Precisely, we consider the following structural setting.
Assume there exist $q_0, q_{1}\in L^{\infty}(\Omega)$, a constant $C_{*}>0$, and real numbers $\alpha_{2}>\alpha_{1}>0$ such that the nonlinear term $q$ satisfies 
\begin{equation}
\begin{aligned}
& \norm{q(\cdot,\abs{z_{1}})z_{1}-q(\cdot,\abs{z_{2}})z_{2} - q_{0}(\cdot)(z_{1}-z_{2}) - q_{1}(\cdot)(\abs{z_{1}}^{\alpha_{1}}z_{1}-\abs{z_{2}}^{\alpha_{1}}z_{2})}_{L^{\infty}(\Omega)} \\ 
& \quad \le C_{*}\left( \abs{z_{1}}^{\alpha_{2}} + \abs{z_{2}}^{\alpha_{2}} \right) \abs{z_{1}-z_{2}}
\end{aligned}  \label{eq:generalized-Kerr}
\end{equation}
\end{subequations}
for all $z_{1},z_{2}\in\mC$ with $\abs{z_{1}}\le 1$ and $\abs{z_{2}}\le 1$. 
Condition \eqref{eq:generalized-Kerr} describes nonlinearities for which
$q_1(x)|z|^{\alpha_1}z$ is the leading nonlinear term near $z=0$,
while the remaining nonlinear part is of strictly higher order.
Roughly speaking, condition \eqref{eq:generalized-Kerr} may be viewed as
\[
q(x,|z|)z
=
q_0(x)z
+
q_1(x)|z|^{\alpha_1}z
+
O(|z|^{1+\alpha_2}),
\qquad z\to0,
\]
where $\alpha_2>\alpha_1$. 
The simplest and typical example of \eqref{eq:generalized-Kerr} is the well-known Kerr-type nonlinearity, which corresponds to the special case 
\begin{equation*}
q(x,\abs{z})=q_{0}(x) + q_{1}(x)\abs{z}^{2}.  
\end{equation*}
In this case, one may take $\alpha_1=2$ and any $\alpha_2>2$, and since the remainder vanishes identically, \eqref{eq:generalized-Kerr} holds for any $C_*>0$. This nonlinearity arises in models of nonlinear optics, see e.g. \cite{Boy20NonlinearOptics,MN19NonlinearOptics} for a physical description. Moreover, this equation arises as the continuum limit of a family of discrete nonlinear Schr{\"o}dinger equations with long-range lattice interactions, modeling charge transport in biopolymers such as DNA, see e.g. \cite{KLS13ContinuumLimitDNLS}. Using \Cref{lem:Kerr-inequality} below, one can readily verify that another example of \eqref{eq:generalized-Kerr} is the Ginzburg-Landau-type nonlinearity, corresponding to the special case
\begin{equation*}
q(x,\abs{z})=q_{0}(x) + q_{1}(x)\abs{z}^{2} + q_{2}(x)\abs{z}^{4}.  
\end{equation*}
Here one may take $\alpha_1=2$, $\alpha_2=4$, and a constant $C_*\geq 16 \|q_2\|_{L^\infty(\Omega)}$. 
This nonlinearity arises in models of light propagation in nonlinear dissipative media, see e.g. \cite{qiu2020soliton}. 

To establish the well-posedness of the forward problem, we additionally assume that
\begin{equation*}
q_{0} \le 0 \quad \text{in $\Omega$.} 
\end{equation*}
Such a sign condition can arise, for example, for soliton envelope functions satisfying the stationary nonlinear Schr{\"o}dinger equation when the soliton parameter is sufficiently large, see, e.g., \cite[(2)]{Flach2005Resonant}. 
Furthermore, the nonlinearity considered in \cite{LLST22fractionalpowerlinearity} is also an instance of the form \eqref{eq:generalized-Kerr}, with \(q_0=0\).

\subsection*{Methodology and key difficulties}
In this paper, we apply the enclosure method to reconstruct the convex hull of the unknown inclusion. The enclosure method was originally introduced by Ikehata in \cite{Ikehata99Enclosure, Ikehata00Enclosure}, and has since been extended to various types of equations and systems (see, e.g., \cite{BHKS2018,BKS2015,I2023,KLS2015,NagaUW2011,UW2008,WZ13Enclosure}).
The central idea of the method is to substitute a series of special solutions into a designed indicator functional; the resulting values --- called the indicator function --- are then analyzed. By observing its limiting behavior, one can detect the location of the unknown inclusion.

The method relies on two main components: the special solutions and the indicator functional. In many classical cases, the special solutions are taken to be complex geometrical optics (CGO) solutions. However, in the previous work \cite{Kuan24Enclosure} involving nonlinear equations with power-type nonlinearities, the CGO solutions were replaced by so-called approximate solutions, which were constructed via a Taylor expansion of the associated solution operator. This construction was inspired by the higher-order linearization framework proposed in \cite{LLLS21powernonlinearity,LLLS21semilinear, KU20semilinear}, and adapted to the setting of the enclosure method.

A main difficulty in our setting is the lack of analyticity (cf. \cite{LLLS21powernonlinearity,LLLS21semilinear, KU20semilinear}) of the nonlinear term $z\mapsto q(x,\abs{z})z$, which renders the Taylor-based construction of approximate solutions in \cite{Kuan24Enclosure} inapplicable. 
We therefore develop an alternative approach that does not rely on a higher-order Taylor expansion of the solution operator.
The estimates required for this approach are based on key inequalities from \cite{GKM22NonlinearFactorizationMethod, KW23FractionalNonlinear}. Further details are provided in \Cref{sec:inverse-problem}. 

\subsection*{Mathematical framework and main theorem}

To formulate the inverse problem in terms of boundary data, it is essential to ensure that the forward problem is well-posed. In particular, the solution to the boundary value problem must exist uniquely and depend continuously on the boundary data, so that the associated Dirichlet-to-Neumann (DtN) map is well-defined. 
A detailed discussion of the well-posedness is provided in \Cref{sec:forward-problem}.

We formulate our problem using the following spaces: 
\begin{equation*}
\mA^{2}(\Omega) := H^{2}(\Omega)\cap L^{\infty}(\Omega), \quad \mA^{\frac{3}{2}}(\partial\Omega) := H^{\frac{3}{2}}(\partial\Omega)\cap L^{\infty}(\partial\Omega)
\end{equation*}
each equipped with the norm 
\begin{equation*}
\norm{\cdot}_{\mA^{2}(\Omega)} := \norm{\cdot}_{H^{2}(\Omega)} + \norm{\cdot}_{L^{\infty}(\Omega)} ,\quad \norm{\cdot}_{\mA^{\frac{3}{2}}(\partial\Omega)} := \norm{\cdot}_{H^{\frac{3}{2}}(\partial\Omega)} + \norm{\cdot}_{L^{\infty}(\partial\Omega)}.
\end{equation*}
Both spaces are Banach spaces under these norms. 
Using \Cref{lem:well-posedness} below, one can identify a small parameter
\begin{equation}
\delta_{0}=\delta_{0}(\Omega,C_{*},\alpha_{1},\alpha_{2},\norm{q_{0}}_{L^{\infty}(\Omega)},\norm{q_{1}}_{L^{\infty}(\Omega)})>0 \label{eq:small-parameter-fixed}
\end{equation}
such that if 
\begin{equation*}
f\in U_{\delta_{0}}:= \left\{ f\in\mA^{\frac{3}{2}}(\partial\Omega) : \norm{f}_{\mA^{\frac{3}{2}}(\partial\Omega)} \le \delta_{0} \right\}, 
\end{equation*}
then there exists a unique ``small'' solution $u\in \mA^{2}(\Omega)$ of \eqref{eq:main}. Since $\partial\Omega\in C^{\infty}$, using the trace theorem \cite[Theorem~7.3.11]{AH09TheoreticalNumericalAnalysis}, we can define the Dirichlet-to-Neumann map (DN-map) as 
\begin{equation}
\Lambda_{q}(f):=\partial_{\nu}u_{f} \in H^{1/2}(\partial\Omega) \quad \text{for all $f\in U_{\delta_{0}}$}, \label{eq:DN-map} 
\end{equation}
where $\nu$ is the unit outward normal vector on $\partial\Omega$ and $\partial_{\nu}=\nu\cdot\nabla$ is the normal derivative. 

Then, with the forward problem properly defined, we move on to the inverse problem studied in this paper: reconstructing the convex hull of an unknown inclusion from boundary measurements in a semilinear setting with non-analytic nonlinearities.
We assume that the unknown inclusion $D$ is an open set satisfying $\overline{D}\subset\Omega$, and represented by a perturbation of the nonlinear coefficient, modeled as
\begin{subequations} \label{eq:perturbation-condition} 
\begin{equation}
q_{1}(x)=q_{1,b}(x) + \chi_{D}q_{1,D}(x) \quad \text{for all $x\in\Omega$,} \label{eq:perturbation} 
\end{equation}
where $q_{1,b},q_{1,D}\in L^{\infty}(\Omega)$, $q_{1,b}$ is a known function (representing the ``perfect'' material), $q_{1,D}$ is unknown (representing a ``defect'' in the material). 

We further assume that the ``defect'' is detectable in the sense, that is we need a jump condition for the inclusion, that 
\begin{equation}
\inf_{x\in D} q_{1,D}(x) \ge \mu \quad \text{or} \quad \sup_{x\in D} q_{1,D}(x) \le -\mu \label{eq:contrast-unknown}
\end{equation}
for some constant $\mu>0$.
\end{subequations} 

Our main contribution is to extend the enclosure-method reconstruction in \cite{Kuan24Enclosure}, developed for power-type nonlinearities, to the generalized Kerr-type setting described by \eqref{eq:generalized-Kerr}. In contrast to the power-type case, the nonlinear term considered here need not possess the higher-order differentiability required for the Taylor-based construction of approximate solutions in \cite{Kuan24Enclosure}. We overcome this difficulty by constructing an approximate solution through a first iteration starting from the solution of the linearized problem. This construction is developed in detail in \Cref{sec:inverse-problem}. A key point is that this approximation already provides sufficient accuracy for the enclosure-method argument. This allows us to reconstruct the convex hull of the inclusion \(D\) from the boundary measurement \(\Lambda_q\), as stated in the following theorem. 

\begin{theorem}[see \Cref{thm:main-detailed} below for a detailed statement] \label{thm:main} 
Let $\Omega \subset \mathbb{R}^n$ ($n \ge 2$) be a bounded $C^{\infty}$-domain and $D$ be an open set satisfying $\overline{D}\subset\Omega$. Assume there exist $0 \ge q_{0}\in H^s(\Omega)$ for some $s>n/2$ with $s>2$, $q_{1}\in L^{\infty}(\Omega)$ satisfying \eqref{eq:perturbation-condition}, a constant $C_{*}>0$, and real numbers $\alpha_{2}>\alpha_{1}>0$, such that the nonlinear term $q$ satisfies \eqref{eq:generalized-Kerr}. Then the convex hull of $D$ is uniquely determined by the Dirichlet-to-Neumann map $\Lambda_{q}(f):U_{\delta_{0}}\rightarrow H^{1/2}(\partial\Omega)$ defined in \eqref{eq:DN-map}. 
\end{theorem}

It should be noted that our method does not require the exact expression of the nonlinearity $q$. This allows for greater flexibility in the admissible forms of $q$.

\subsection*{Organization}
\addtocontents{toc}{\SkipTocEntry}

In \Cref{sec:forward-problem}, we establish the well-posedness of the boundary value problem \eqref{eq:main}, followed by the proof of our main result (\Cref{thm:main}) in \Cref{sec:inverse-problem}. Finally, we address the numerical challenges associated with the practical implementation of the enclosure method in \Cref{sec:numerical}.

\section{Forward problem\label{sec:forward-problem}}

Before we tackle the inverse problem, we first consider the forward problem by establishing the well-posedness of the boundary value problem \eqref{eq:main}. We now introduce the following assumption, see also \cite[Assumption~2.1]{GKM22NonlinearFactorizationMethod}:

\begin{assumption}[{\cite[Assumption~1.1]{KW23FractionalNonlinear}}]\label{assu:generalized-Kerr}
There exists $C_{\rm Kerr}>0$, $0 \ge q_{0}\in L^{\infty}(\Omega)$ and a $\alpha>0$ such that 
\begin{equation*}
\norm{q(\cdot,\abs{z_{1}})z_{1}-q(\cdot,\abs{z_{2}})z_{2}-q_{0}(\cdot)(z_{1}-z_{2})}_{L^{\infty}(\Omega)}\le C_{\rm Kerr} \left( \abs{z_{1}}^{\alpha} + \abs{z_{2}}^{\alpha} \right) \abs{z_{1}-z_{2}}
\end{equation*}
for all $z_{1},z_{2}\in\mC$ with $\abs{z_{1}}\le 1$ and $\abs{z_{2}}\le 1$. 
\end{assumption}

We note that condition \eqref{eq:generalized-Kerr} specifies the leading nonlinear term
$q_1(x)|z|^{\alpha_1}z$ and requires the remaining part to be of higher order,
since $\alpha_2>\alpha_1$. The following lemma shows that
\eqref{eq:generalized-Kerr} with $q_0\leq 0$ implies \Cref{assu:generalized-Kerr}
with $\alpha=\alpha_1$. In particular, one may take
$C_{\mathrm{Kerr}}
=C_*+2^{\max\{1,\alpha_1\}}\norm{q_1}_{L^\infty(\Omega)}$.

\begin{lemma}[{\cite[Lemma~A.1]{GKM22NonlinearFactorizationMethod}}]\label{lem:Kerr-inequality}
The inequality 
\begin{equation*}
\abs{\abs{a}^{\alpha}a-\abs{b}^{\alpha}b} \le 2(\abs{a}+\abs{b})^{\alpha} \abs{a-b}
\end{equation*}
holds for all $a,b\in\mC$ and $\alpha>0$. 
\end{lemma}

Since $0 \ge q_{0}\in L^{\infty}(\Omega)$, for each $f\in \mA^{\frac{3}{2}}(\partial\Omega)$ one can use Fredholm theory to guarantee that there exists a unique solution $v_{f} \in H^{2}(\Omega)$ to 
\begin{equation}
(\Delta + q_{0}) v_{f} = 0 \text{ in $\Omega$} ,\quad \left.v_{f}\right|_{\partial\Omega} = f. \label{eq:linear-part}
\end{equation} 

Let $\vartheta\in\mR$ be a number, and we see that 
\begin{equation*}
(\Delta + q_{0}) \Re(e^{\bfi\vartheta}v_{f}) = 0 \text{ in $\Omega$} ,\quad \left.\Re(e^{\bfi\vartheta}v_{f})\right|_{\partial\Omega} = \Re(e^{\bfi\vartheta}f). 
\end{equation*} 
Define $v_{*}:=\norm{f}_{L^{\infty}(\partial\Omega)}$. We now see that $v_{*}\pm \Re(e^{\bfi\vartheta}v_{f}) \in H^{1}(\Omega)$ solves 
\begin{equation*}
-(\Delta+q_{0})(v_{*}\pm \Re(e^{\bfi\vartheta}v_{f})) = -q_{0}v_{*} \ge0 \text{ in $\Omega$.} 
\end{equation*}
Since $(v_{*}\pm \Re(e^{\bfi\vartheta}v_{f}))|_{\partial\Omega}=\norm{f}_{L^{\infty}(\partial\Omega)}\pm\Re(e^{\bfi\vartheta}f) \ge 0$ and $q_{0}\le 0$ a.e. in $\Omega$, using maximum principle we see that $v_{*}\pm \Re(e^{\bfi\vartheta}v_{f})\ge 0$ a.e. in $\Omega$, that is, 
\begin{equation*}
\abs{\Re(e^{\bfi\vartheta}v_{f})(x)} \le \norm{f}_{L^{\infty}(\partial\Omega)} \quad \text{for all $\vartheta\in\mR$ and a.e. $x\in\Omega$.} 
\end{equation*}
For a.e. $x\in\Omega$, choosing $\vartheta=-\arg v_{f}(x)$ yields  
\begin{equation}
\norm{v_{f}}_{L^{\infty}(\Omega)} \le \norm{f}_{L^{\infty}(\partial\Omega)}. \label{eq:Linfty-u-linear-part}
\end{equation} 
By writing \eqref{eq:linear-part} as 
\begin{equation*}
\Delta v_{f} = -q_{0}v_{f} \text{ in $\Omega$} ,\quad \left.v_{f}\right|_{\partial\Omega} = f, 
\end{equation*}
the trace theorem \cite[Theorem~7.3.11]{AH09TheoreticalNumericalAnalysis} and the elliptic estimate \cite[Theorem~8.12]{GT01Elliptic} together imply 
\begin{equation}
\norm{v_{f}}_{H^{2}(\Omega)} \le C\left(\norm{f}_{H^{\frac{3}{2}}(\partial\Omega)} + (1+\norm{q_{0}}_{L^{\infty}(\Omega)})\norm{v_{f}}_{L^{2}(\Omega)}\right). \label{eq:H2-u-linear-part}
\end{equation}
Combining \eqref{eq:H2-u-linear-part} and \eqref{eq:Linfty-u-linear-part}, we conclude 
\begin{equation}
\norm{v_{f}}_{\mA^{2}(\Omega)} \le C (1 + \norm{q_{0}}_{L^{\infty}(\Omega)})\norm{f}_{\mA^{\frac{3}{2}}(\partial\Omega)} \label{eq:Linfty-u-linear-part-estimate}
\end{equation}
for some constant $C=C(\Omega)>1$ which is independent of $q_{0}$.

For each $F \in L^{2}(\Omega)$, again one can use the Lax-Milgram theorem to guarantee that there exists a unique solution $w_{F} \in H_{0}^{1}(\Omega)$ of 
\begin{equation}
(\Delta + q_{0})w_{F} = F \text{ in $\Omega$} ,\quad \left.w_{F}\right|_{\partial\Omega}=0. \label{eq:solution-operator-contraction1}
\end{equation}
Since $q_{0}\le 0$, the Poincar{\'e} inequality implies  
\begin{equation*} 
\scalebox{0.9}{$\begin{aligned} 
& \norm{\Re w_{F}}_{H_{0}^{1}(\Omega)}^{2} \le C \norm{\nabla\Re w_{F}}_{L^{2}(\Omega)}^{2} = C\left(\int_{\Omega} q_{0}\abs{\Re w_{F}}^{2} \, \rmd x - \int_{\Omega} (\Re F)(\Re w_{F}) \, \rmd x\right) \le - C\int_{\Omega} (\Re F)(\Re w_{F}) \, \rmd x, \\ 
& \norm{\Im w_{F}}_{H_{0}^{1}(\Omega)}^{2} \le C \norm{\nabla\Im w_{F}}_{L^{2}(\Omega)}^{2} = C\left(\int_{\Omega} q_{0}\abs{\Im w_{F}}^{2} \, \rmd x - \int_{\Omega} (\Im F)(\Im w_{F}) \, \rmd x\right) \le - C\int_{\Omega} (\Im F)(\Im w_{F}) \, \rmd x, 
\end{aligned}$} 
\end{equation*}
and thus we conclude that 
\begin{equation}
\begin{aligned} 
& \norm{w_{F}}_{H_{0}^{1}(\Omega)} = \left( \norm{\Re w_{F}}_{H_{0}^{1}(\Omega)}^{2} + \norm{\Im w_{F}}_{H_{0}^{1}(\Omega)}^{2} \right)^{1/2} \\ 
&\quad \le \left( C\norm{\Re F}_{L^{2}(\Omega)}^{2} + C\norm{\Im F}_{L^{2}(\Omega)}^{2} \right)^{1/2} = C\norm{F}_{L^{2}(\Omega)} 
\end{aligned} \label{eq:H1-estimate} 
\end{equation} 
for some constant $C=C(\Omega)>1$ which is independent of $q_{0}$. We now write \eqref{eq:solution-operator-contraction1} as 
\begin{equation*}
\Delta w_{F} = F - q_{0}w_{F} \text{ in $\Omega$} ,\quad \left.w_{F}\right|_{\partial\Omega}=0, 
\end{equation*}
and we use the elliptic estimate \cite[Theorem~8.12]{GT01Elliptic} to see that 
\begin{equation*}
\norm{w_{F}}_{H^{2}(\Omega)} \le C \left( \norm{F}_{L^{2}(\Omega)} + (1+\norm{q_{0}}_{L^{\infty}(\Omega)})\norm{w_{F}}_{L^{2}(\Omega)} \right),
\end{equation*}
hence from \eqref{eq:H1-estimate} we reach 
\begin{equation}
\norm{w_{F}}_{H^{2}(\Omega)} \le C(1+\norm{q_{0}}_{L^{\infty}(\Omega)})\norm{F}_{L^{2}(\Omega)} \label{eq:H2-estimate} 
\end{equation}
for some constant $C=C(\Omega)>1$ which is independent of $q_{0}$. 
We now further assume that $F\in L^{\infty}(\Omega)$ and let $w_{*}^{\Re} \in H_{0}^{1}(\Omega)$ be the unique solution to 
\begin{equation*}
-\Delta w_{*}^{\Re} = \norm{\Re F}_{L^{\infty}(\Omega)} \text{ in $\Omega$} ,\quad \left. w_{*}^{\Re}\right|_{\partial\Omega}=0. 
\end{equation*}
By using the maximum principle for weak solutions \cite[Theorem~8.19]{GT01Elliptic} one sees that 
\begin{equation}
w_{*}^{\Re} \ge 0 \quad \text{a.e. in $\Omega$.} \label{eq:non-negative-v-star}
\end{equation}
In addition, by using \cite[Theorem~8.16]{GT01Elliptic}, there exists a constant $C=C(\Omega)>1$ such that 
\begin{equation}
\norm{w_{*}^{\Re}}_{L^{\infty}(\Omega)} \le C\norm{\Re F}_{L^{\infty}(\Omega)}. \label{eq:Linfty-bound1}
\end{equation}
We now see that $w_{*}^{\Re}\pm\Re w_{F} \in H_{0}^{1}(\Omega)$ is the unique solution to 
\begin{equation*}
-(\Delta + q_{0})(w_{*}^{\Re}\pm\Re w_{F}) = \norm{\Re F}_{L^{\infty}(\Omega)} \mp\Re F - q_{0}w_{*}^{\Re} \text{ in $\Omega$} ,\quad \left.\left(w_{*}^{\Re}\pm\Re w_{F}\right)\right|_{\partial\Omega} = 0. 
\end{equation*}
Since $q_{0}\le 0$, then from \eqref{eq:non-negative-v-star} we have $\norm{\Re F}_{L^{\infty}(\Omega)} \mp\Re F - q_{0}w_{*}^{\Re} \ge 0$ a.e. in $\Omega$, then by using the maximum principle for weak solutions \cite[Theorem~8.19]{GT01Elliptic} and \eqref{eq:non-negative-v-star} one sees that $w_{*}^{\Re}\pm\Re w_{F} \ge 0$ a.e. in $\Omega$, that is, $\abs{\Re w_{F}} \le w_{*}^{\Re}=\abs{w_{*}^{\Re}}$ a.e. in $\Omega$, together with \eqref{eq:Linfty-bound1} we now reach 
\begin{equation*}
\norm{\Re w_{F}}_{L^{\infty}(\Omega)} \le C\norm{\Re F}_{L^{\infty}(\Omega)}. 
\end{equation*}
An analogous argument yields 
\begin{equation*}
\norm{\Im w_{F}}_{L^{\infty}(\Omega)} \le C\norm{\Im F}_{L^{\infty}(\Omega)}. 
\end{equation*}
Therefore, combining the above two estimates, we obtain 
\begin{equation} 
\norm{w_{F}}_{L^{\infty}(\Omega)} \le \sqrt{2}C\norm{F}_{L^{\infty}(\Omega)}. \label{eq:Linfty-bound2}
\end{equation}
Now from \eqref{eq:H2-estimate} and \eqref{eq:Linfty-bound2} we can define the bounded linear operator 
\begin{equation}
\calS : L^{\infty}(\Omega) \rightarrow H_{0}^{1}(\Omega)\cap \mA^{2}(\Omega) ,\quad \calS[F]:=w_{F}, 
\end{equation}
where $w_{F}$ is the unique solution to \eqref{eq:solution-operator-contraction1} satisfying 
\begin{equation}
\norm{\calS[F]}_{\mA^{2}(\Omega)} \le C(1+\norm{q_{0}}_{L^{\infty}(\Omega)}) \norm{F}_{L^{\infty}(\Omega)} \label{eq:contraction-operator-S}
\end{equation}
for some constant $C=C(\Omega)>1$ which is independent of $q_{0}$.

For each $\delta>0$, let us define the closed ball 
\begin{equation*}
X_{\delta} := \left\{ w\in H_{0}^{1}(\Omega) \cap \mA^{2}(\Omega) : \norm{w}_{\mA^{2}(\Omega)} \le \delta \right\}. 
\end{equation*}
We are now ready to prove the well-posedness of \eqref{eq:generalized} under \Cref{assu:generalized-Kerr} following the ideas in \cite[Theorem~2.1]{KW23FractionalNonlinear}, which is based on the contraction mapping theorem. 

\begin{lemma}\label{lem:well-posedness}
Let $0\ge q_{0} \in L^{\infty}(\Omega)$ and let $q$ satisfy \Cref{assu:generalized-Kerr}. There exists a sufficiently small $\delta_{0}=\delta_{0}(\alpha,C_{\rm Kerr},\Omega,\norm{q_{0}}_{L^{\infty}(\Omega)})\in(0,1)$ such that the following statement holds true: if $f\in U_{\delta_{0}}$, then there exists a unique solution $u_{f}\in \mA^{2}(\Omega)$ to \eqref{eq:generalized} satisfying $u_{f}-v_{f}\in X_{\delta}$ with $\delta=\norm{f}_{\mA^{\frac{3}{2}}(\partial\Omega)}$ and 
\begin{subequations}\label{eq:apriori1}
\begin{equation}
\norm{u_{f}-v_{f}}_{\mA^{2}(\Omega)} \le C_{\rm Kerr} (2C^{2})^{1+\alpha}(1+\norm{q_{0}}_{L^{\infty}(\Omega)})^{2+\alpha}\norm{f}_{\mA^{\frac{3}{2}}(\partial\Omega)}^{1+\alpha} \label{eq:apriori-solution-u-1}
\end{equation}
and 
\begin{equation}
\norm{u_{f}}_{\mA^{2}(\Omega)} \le C (1 + \norm{q_{0}}_{L^{\infty}(\Omega)})\norm{f}_{\mA^{\frac{3}{2}}(\partial\Omega)} \label{eq:apriori-solution-u-2}
\end{equation}
\end{subequations} 
for some constant $C=C(\Omega)>1$ which is independent of $q_{0},\alpha,C_{\rm Kerr}$, where $v_{f}$ is the unique solution to \eqref{eq:linear-part}. 
\end{lemma}

\begin{proof}
From \eqref{eq:Linfty-u-linear-part-estimate} we have 
\begin{equation}
\norm{v_{f}}_{\mA^{2}(\Omega)} \le C (1 + \norm{q_{0}}_{L^{\infty}(\Omega)}) \delta. \label{eq:Linfty-u-linear-part-estimate-delta}
\end{equation}
We define 
\begin{equation*}
\calF[w] := -\left(q(\cdot,\abs{w+v_{f}})-q_{0}\right)(w+v_{f}). 
\end{equation*}
We first show that 
\begin{equation}
\calS [\calF[w]] \in X_{\delta} \quad \text{for all $w\in X_{\delta}$}. \label{eq:contraction1} 
\end{equation}
For each $w\in X_{\delta}$, we have 
\begin{equation}
\norm{w+v_{f}}_{L^{\infty}(\Omega)}\le\norm{w}_{L^{\infty}(\Omega)}+\norm{v_{f}}_{L^{\infty}(\Omega)}\le \delta + C(1+\norm{q_{0}}_{L^{\infty}(\Omega)})\delta \le 1 \label{eq:verification}
\end{equation}
provided that $\delta_{0}>0$ is chosen sufficiently small. Thus, this choice $z_{1}=w+v_{f}$ and $z_{2}=0$ is admissible for \Cref{assu:generalized-Kerr}. Applying \eqref{eq:contraction-operator-S}, we obtain 
\begin{equation}
\begin{aligned}
& \norm{\calS [\calF[w]]}_{\mA^{2}(\Omega)} \le C(1+\norm{q_{0}}_{L^{\infty}(\Omega)}) \norm{\calF[w]}_{L^{\infty}(\Omega)} \\ 
& \quad \le C_{\rm Kerr} C(1+\norm{q_{0}}_{L^{\infty}(\Omega)}) \norm{w+v_{f}}_{L^{\infty}(\Omega)}^{1+\alpha} \\
& \quad \le C_{\rm Kerr} C^{2+\alpha}2^{1+\alpha}(1+\norm{q_{0}}_{L^{\infty}(\Omega)})^{2+\alpha}\delta^{1+\alpha}
\end{aligned} \label{eq:apriori-estimate1}
\end{equation}
where the last inequality follows from \eqref{eq:Linfty-u-linear-part-estimate-delta}. If we choose $\delta_{0}>0$ such that 
\begin{equation}
C_{\rm Kerr} C^{2+\alpha}2^{1+\alpha}(1+\norm{q_{0}}_{L^{\infty}(\Omega)})^{2+\alpha}\delta_{0}^{\alpha} \le 1, \label{eq:small-delta-0-condition1}
\end{equation}
then we conclude \eqref{eq:contraction1}. Next, we want to show that 
\begin{equation}
\text{$\calS[\calF[\cdot]]$ is a contraction map on $X_{\delta}$.} \label{eq:contraction2} 
\end{equation}
Let $w_{1},w_{2}\in X_{\delta}$ and compute 
\begin{equation*}
\begin{aligned}
\calF[w_{1}]-\calF[w_{2}] &= q(\cdot,\abs{w_{2}+v_{f}})(w_{2}+v_{f}) - q(\cdot,\abs{w_{1}+v_{f}})(w_{1}+v_{f}) \\ 
& \qquad - q_{0}((w_{2}+v_{f})-(w_{1}+v_{f})). 
\end{aligned}
\end{equation*}
Since $\calS$ is linear, applying \eqref{eq:contraction-operator-S} and then \Cref{assu:generalized-Kerr} with the choice $z_{1}=w_{2}+v_{f}$ and $z_{2}=w_{1}+v_{f}$ (whose admissibility follows from an argument analogous to that in \eqref{eq:verification}), we obtain 
\begin{equation*}
\begin{aligned}
& \norm{\calS[\calF[w_{1}]]-\calS[\calF[w_{2}]]}_{\mA^{2}(\Omega)} \\
& \quad \le C(1+\norm{q_{0}}_{L^{\infty}(\Omega)}) \norm{\calF[w_{1}]-\calF[w_{2}]}_{L^{\infty}(\Omega)} \\
& \quad \le C_{\rm Kerr} C(1+\norm{q_{0}}_{L^{\infty}(\Omega)}) \left( \norm{w_{1}+v_{f}}_{L^{\infty}(\Omega)}^{\alpha} + \norm{w_{2}+v_{f}}_{L^{\infty}(\Omega)}^{\alpha} \right) \norm{w_{1}-w_{2}}_{L^{\infty}(\Omega)} \\ 
& \quad \le \overbrace{C_{\rm Kerr} C^{1+\alpha}2^{1+\alpha}(1+\norm{q_{0}}_{L^{\infty}(\Omega)})^{1+\alpha}\delta^{\alpha}}^{<1 \text{ since \eqref{eq:small-delta-0-condition1} and $C>1$}} \norm{w_{1}-w_{2}}_{L^{\infty}(\Omega)}, 
\end{aligned}
\end{equation*}
then we conclude \eqref{eq:contraction2} from \eqref{eq:small-delta-0-condition1}. 

From \eqref{eq:contraction1} and \eqref{eq:contraction2}, the Banach fixed point theorem guarantees that there exists a unique $w_{f}\in X_{\delta}$ such that 
\begin{equation}
w_{f} = \calS[\calF[w_{f}]]. \label{eq:fixed-point-contraction}
\end{equation}
We see that $w_{f}$ satisfies 
\begin{equation*}
(\Delta+q_{0})w_{f}=\calF[w_{f}]=-\left(q(\cdot,\abs{w_{f}+v_{f}})-q_{0}\right)(w_{f}+v_{f}) \text{ in $\Omega$} ,\quad w_{f}|_{\partial\Omega} = 0. 
\end{equation*}
We now see that the function $u := w_{f}+v_{f} \in \mA^{2}(\Omega)$ is the unique solution to \eqref{eq:generalized} with $u-v_{f}\in X_{\delta}$. 

On the other hand, from \eqref{eq:apriori-estimate1} and \eqref{eq:fixed-point-contraction} we conclude \eqref{eq:apriori-solution-u-1}. From \eqref{eq:apriori-estimate1} and \eqref{eq:small-delta-0-condition1} as well as \eqref{eq:fixed-point-contraction}, we see that 
\begin{equation*}
\norm{w_{f}}_{\mA^{2}(\Omega)} \le \delta = \norm{f}_{\mA^{\frac{3}{2}}(\partial\Omega)}. 
\end{equation*}
Combining the above inequality with \eqref{eq:Linfty-u-linear-part-estimate}, we conclude \eqref{eq:apriori-solution-u-2}. 
\end{proof}

\section{Proof of the main result\label{sec:inverse-problem}} 

In this section, we adapt the enclosure-method framework to the present non-analytic setting for the reconstruction of the convex hull of the inclusion $D$. The precise detection criterion is established in \Cref{thm:main-detailed}, from which \Cref{thm:main} follows.

A key observation in \cite{Kuan24Enclosure} is that, in order to implement the enclosure method for a nonlinear equation, one does not necessarily need to construct exact CGO solutions to the target nonlinear equation. Instead, one may work with suitably chosen approximate solutions, together with an auxiliary indicator functional, provided that the corresponding approximation errors are sufficiently small for the asymptotic analysis. In \cite{Kuan24Enclosure}, such approximate solutions were constructed through a Taylor expansion of the solution to the nonlinear boundary value problem with respect to the boundary data. This Taylor-based construction is closely related to the higher-order linearization method and relies on the higher-order differentiability of the associated solution map.

However, the situation considered here is different. For generalized Kerr-type nonlinearities satisfying \eqref{eq:generalized-Kerr}, the nonlinear term does not need to possess the higher-order differentiability required for Taylor-based construction in \cite{Kuan24Enclosure}. We therefore seek a different approximation mechanism. Our aim is not to obtain a full higher-order expansion of the nonlinear solution, but rather to construct a sufficiently accurate linear problem approximating the target equation so that the enclosure method can still be carried out. The core idea is to perform the first iteration starting from the linearized solution $v_f$, defined in \eqref{eq:linear-part}. More precisely, we replace $u_f$ in the nonlinear term of \eqref{eq:generalized} by $v_f$, which leads to an approximating linear problem. We denote the solution of this problem by $\widetilde{u}_f$ and refer to it as an approximate solution.

The point of this construction is that this first iteration already provides a sufficiently accurate approximation for the subsequent implementation of the enclosure method. In the following subsections, we develop this idea step by step and show how it leads to the desired enclosure-type reconstruction. Since our construction relies on the well-posedness of the forward problem, we retain the additional assumption \(q_0\leq 0\).

\subsection{Approximate solutions\label{sec:approximate-solution}}
\addtocontents{toc}{\SkipTocEntry}

Let 
\[
\tilde{q}(x,\abs{z}) = q(x,\abs{z}) - q_0(x), \quad \forall x \in \Omega \mbox{ and } z \in \mathbb{C}.
\] 
Recall that $v_f$ denotes the unique solution of \eqref{eq:linear-part}. 
We define $\tilde{u}_f$ as the unique solution to the following \emph{linear} boundary value problem: 
\begin{equation}
\Delta\tilde{u}_{f} + q_{0}(x)\tilde{u}_{f} = -\tilde{q}(x,\abs{v_{f}})v_{f} \text{ in } \Omega ,\quad \left.\tilde{u}_{f}\right|_{\partial\Omega}=f. \label{eq:approximate-solution}
\end{equation} 
We refer to $\widetilde{u}_f$ as the approximate solution associated with $f$. 
The following lemma provides an estimate for the approximation error between $\widetilde{u}_{f}$ and $u_{f}$, which will be used in the subsequent construction of the indicator functional.

\begin{lemma}\label{lem:linearization-solution-operator} 
Let $q$ satisfy \eqref{eq:generalized-Kerr}. Then there exists a constant $C=C(\Omega,\norm{q_{0}}_{L^{\infty}(\Omega)},\norm{q_{1}}_{L^{\infty}(\Omega)},\alpha_1,\alpha_2, C_*)>0$ such that 
\begin{equation*}
\norm{u_{f}-\tilde{u}_{f}}_{\mA^{2}(\Omega)} \le C\norm{f}_{\mA^{\frac{3}{2}}(\partial\Omega)}^{1+2\alpha_{1}} \quad \text{for all $f\in U_{\delta_{0}}$,}
\end{equation*}
where $\delta_{0}\in(0,1)$ denotes the constant specified in \Cref{lem:well-posedness}. 
\end{lemma}

\begin{proof} 
Since $u_{f}-\tilde{u}_{f}=-\calS\left[ \tilde{q}(\cdot,\abs{u_{f}})u_{f}-\tilde{q}(\cdot,\abs{v_{f}})v_{f} \right]$, in other words, 
\begin{equation*}
(\Delta+q_{0})(u_{f}-\tilde{u}_{f}) = -\left(\tilde{q}(x,\abs{u_{f}})u_{f}-\tilde{q}(x,\abs{v_{f}})v_{f}\right) \text{ in $\Omega$} ,\quad \left.(u_{f}-\tilde{u}_{f})\right|_{\partial\Omega}=0. 
\end{equation*}
Using \eqref{eq:contraction-operator-S}, we then obtain that 
\begin{equation*}
\begin{aligned}
& \norm{u_{f}-\tilde{u}_{f}}_{\mA^{2}(\Omega)} \le C(1+\norm{q_{0}}_{L^{\infty}(\Omega)}) \norm{\tilde{q}(\cdot,\abs{u_{f}})u_{f}-\tilde{q}(\cdot,\abs{v_{f}})v_{f}}_{L^{\infty}(\Omega)} \\ 
& \quad \le C(1+\norm{q_{0}}_{L^{\infty}(\Omega)}) \left(\begin{aligned}
& \norm{\tilde{q}(\cdot,\abs{u_{f}})u_{f}-\tilde{q}(\cdot,\abs{v_{f}})v_{f}-q_{1}(\cdot)(\abs{u_{f}}^{\alpha_{1}}u_{f}-\abs{v_{f}}^{\alpha_{1}}v_{f})}_{L^{\infty}(\Omega)} \\ 
& +\norm{q_{1}(\cdot)(\abs{u_{f}}^{\alpha_{1}}u_{f}-\abs{v_{f}}^{\alpha_{1}}v_{f})}_{L^{\infty}(\Omega)}
\end{aligned}\right). 
\end{aligned}
\end{equation*}
From \eqref{eq:apriori-solution-u-2} and \eqref{eq:Linfty-u-linear-part-estimate-delta}, the smallness of $\delta_{0}\in(0,1)$ ensures that 
\begin{equation*}
\norm{u_{f}}_{L^{\infty}(\Omega)} \le 1 \quad \text{and} \quad \norm{v_{f}}_{L^{\infty}(\Omega)} \le 1
\end{equation*}
which allows us to apply \eqref{eq:generalized-Kerr}. 
Consequently, we use \eqref{eq:generalized-Kerr} and the stability estimates in \eqref{eq:apriori1} and \eqref{eq:Linfty-u-linear-part-estimate-delta}, to see that 
\begin{equation*}
\begin{aligned}
& \norm{u_{f}-\tilde{u}_{f}}_{\mA^{2}(\Omega)} \le C(1+\norm{q_{0}}_{L^{\infty}(\Omega)})\left(\begin{aligned}
& \norm{(\abs{u_{f}}^{\alpha_{2}}+\abs{v_{f}}^{\alpha_{2}})(u_{f}-v_{f})}_{L^{\infty}(\Omega)} \\ 
& +2^{\max\{1,\alpha_{1}\}} \norm{q_{1}}_{L^{\infty}(\Omega)}\norm{(\abs{u_{f}}^{\alpha_{1}}+\abs{v_{f}}^{\alpha_{1}})(u_{f}-v_{f})}_{L^{\infty}(\Omega)}
\end{aligned}\right) \\ 
& \quad \le C(1+\norm{q_{0}}_{L^{\infty}(\Omega)})^{3+\alpha_1}\left(\begin{aligned}
& (1+\norm{q_{0}}_{L^{\infty}(\Omega)})^{\alpha_{2}}\norm{f}_{\mA^{\frac{3}{2}}(\partial\Omega)}^{\alpha_{2}} \\ 
& +(1+\norm{q_{0}}_{L^{\infty}(\Omega)})^{\alpha_{1}}\norm{q_{1}}_{L^{\infty}(\Omega)}\norm{f}_{\mA^{\frac{3}{2}}(\partial\Omega)}^{\alpha_{1}}
\end{aligned}\right)\norm{f}_{\mA^{\frac{3}{2}}(\partial\Omega)}^{1+\alpha_{1}} \\ 
& \quad \le C(\Omega,\norm{q_{0}}_{L^{\infty}(\Omega)},\alpha_1,\alpha_2,\norm{q_{1}}_{L^{\infty}(\Omega)}) \norm{f}_{\mA^{\frac{3}{2}}(\partial\Omega)}^{1+2\alpha_{1}}, 
\end{aligned}
\end{equation*}
which conclude our lemma. 
\end{proof}

\subsection{An indicator functional and its approximation}  
\addtocontents{toc}{\SkipTocEntry}

As in the classical enclosure method, the reconstruction is based on comparing the boundary response of the medium containing the inclusion with that of a known background medium. 
Recalling that $q_{1,b}$ is the known background coefficient introduced in \eqref{eq:perturbation}, we define a reference background model by
\[
q_{b}(x,\abs{z}):=q_{0}(x) + q_{1,b}\abs{z}^{\alpha_{1}},\quad \forall x\in\Omega \mbox{ and } z\in\mC,
\]
and consider the corresponding background problem 
\begin{equation}\label{eq:background}
\Delta u_{b,f} + q_{b}(x,\abs{u_{b,f}})u_{b,f} = 0 \text{ in $\Omega$} ,\quad u_{b,f}|_{\partial\Omega}=f.
\end{equation}
Applying \Cref{lem:well-posedness} to the background problem \eqref{eq:background}, and reducing $\delta_{0}$ if necessary, we see that for every $f\in U_{\delta_{0}}$, there exists a unique solution $u_{b,f}\in\mA^{2}(\Omega)$ to \eqref{eq:background}.

Since the background model has the same linearized operator $\Delta+q_0$ as the target equation, the solution $v_f$ to \eqref{eq:linear-part} is also the linearized solution of the background problem \eqref{eq:background}.
We then apply the same approximate-solution construction as in \Cref{sec:approximate-solution} to this background problem. 
Using the same construction, let $\widetilde{u}_{b,f}$ be the unique solution to
\begin{equation}
\Delta \widetilde{u}_{b,f}
+ q_{0}(x)\widetilde{u}_{b,f}
=
-q_{1,b}(x)\abs{v_{f}}^{\alpha_{1}}v_{f}
\quad \text{in } \Omega,
\qquad
\left.\widetilde{u}_{b,f}\right|_{\partial\Omega}=f. \label{eq:background-linearized}
\end{equation} 
Applying \Cref{lem:linearization-solution-operator} to the background problem, there exists a constant
\[
C=C\bigl(
\Omega,
\norm{q_{0}}_{L^{\infty}(\Omega)},
\norm{q_{1,b}}_{L^{\infty}(\Omega)},
\alpha_{1}
\bigr)>0
\]
such that
\begin{equation*}
\norm{u_{b,f}-\widetilde{u}_{b,f}}_{\mA^{2}(\Omega)}
\le
C\norm{f}_{\mA^{\frac{3}{2}}(\partial\Omega)}^{1+2\alpha_{1}}.
\end{equation*}

Having constructed the approximate solution for the reference background problem, we now introduce the corresponding indicator functional.
Motivated by Ikehata's enclosure method \cite{Ikehata99Enclosure,Ikehata00Enclosure} and the approximate-solution approach in \cite{Kuan24Enclosure}, we define 
\begin{equation*}
\mI(f) := \int_{\partial\Omega} (\partial_{\nu}u_{f}-\partial_{\nu}\tilde{u}_{b,f})\overline{f}\,\rmd S \quad \text{for all $f\in U_{\delta_{0}}$}. 
\end{equation*}
As in \cite{Kuan24Enclosure}, we further introduce the auxiliary approximate indicator functional 
\begin{equation*}
\widetilde{\mI}(f) := \int_{\partial\Omega} (\partial_{\nu}\tilde{u}_{f}-\partial_{\nu}\tilde{u}_{b,f})\overline{f}\,\rmd S \quad \text{for all $f\in U_{\delta_{0}}$.} 
\end{equation*}
Notice that the expression in $\widetilde{\mI}$ involves only solutions of linear equations. 
The following lemma provides an estimate for the error in approximating $\mI$ by $\widetilde{\mI}$.

\begin{lemma}\label{lem:linearization-indicator}
Let $q$ satisfy \eqref{eq:generalized-Kerr}. Then there exists a constant $C=C(\Omega,\norm{q_{0}}_{L^{\infty}(\Omega)},\norm{q_{1}}_{L^{\infty}(\Omega)},\alpha_1,C_*)>0$ such that
\begin{equation*}
\abs{\mI(f)-\widetilde{\mI}(f)} \le C\norm{f}_{\mA^{\frac{3}{2}}(\partial\Omega)}^{2+2\alpha_{1}} \quad \text{for all $f\in U_{\delta_{0}}$,} 
\end{equation*}
where $\delta_{0}\in(0,1)$ is is as specified above.
\end{lemma}

\begin{proof}
Using the trace theorem \cite[Theorem~7.3.11]{AH09TheoreticalNumericalAnalysis} (which requires condition $\partial\Omega\in C^{1,1}$), for each $f\in U_{\delta_{0}}$, we compute that 
\begin{equation*}
\begin{aligned}
& \abs{\mI(f)-\widetilde{\mI}(f)} = \left| \int_{\partial\Omega} (\partial_{\nu}u_{f}-\partial_{\nu}\tilde{u}_{f})\overline{f}\,\rmd S \right| \\
& \quad \le \norm{\partial_{\nu}u_{f}-\partial_{\nu}\tilde{u}_{f}}_{L^{2}(\partial\Omega)}\norm{f}_{L^{2}(\partial\Omega)} \\
& \quad \le C\norm{u_{f}-\tilde{u}_{f}}_{H^{2}(\Omega)}\norm{f}_{L^{2}(\partial\Omega)}. 
\end{aligned}
\end{equation*}
Our lemma follows directly from \Cref{lem:linearization-solution-operator}. 
\end{proof}

We now prove that $\widetilde{\mI}$ can be expressed solely in terms of the unique solution $v_{f}$ to \eqref{eq:linear-part}.

\begin{lemma}\label{lem:simplification}
Let $q$ satisfy \eqref{eq:generalized-Kerr}. 
Define 
\begin{equation*}
R(x,\abs{z}):=q(x,\abs{z})-q_{0}(x)-q_{1}(x)\abs{z}^{\alpha_{1}} \quad \text{for all $x\in\Omega$ and for all $z\in\mC$,} 
\end{equation*}
then $\widetilde{\mI}$ can be written as
\begin{equation*}
\widetilde{\mI}(f) = -\int_{D} q_{1,D}(x)\abs{v_{f}}^{2+\alpha_{1}}\,\rmd x - \int_{\Omega} R(x,\abs{v_{f}})\abs{v_{f}}^{2}\,\rmd x \quad \text{for all $f\in U_{\delta_{0}}$,}
\end{equation*}
where $v_{f}$ is the unique solution to \eqref{eq:linear-part} and $\delta_{0}\in(0,1)$ is as specified above.
\end{lemma}

\begin{remark*}
By choosing $z_{1}=v_{f}(x)$ and $z_{2}=0$ in \eqref{eq:generalized-Kerr}, one sees that 
\begin{equation}
\abs{R(x,\abs{v_{f}(x)})} \le C_{*}\abs{v_{f}}^{\alpha_{2}}. \label{eq:remainder-estimate-2} 
\end{equation}
Since $\alpha_2>\alpha_1$, the second term in the representation of $\widetilde{\mI}(f)$ is of higher order in $v_f$ than the first term. 
Thus, when $v_f$ is sufficiently small, the first term provides the leading contribution, whereas the term involving $R$ can be treated as a remainder. 
\end{remark*}

\begin{proof}[Proof of \Cref{lem:simplification}]
We see that the function $z_{f}:=\tilde{u}_{f}-\tilde{u}_{b,f}$ satisfies 
\begin{equation*}
(\Delta+q_{0})z_{f} = -\chi_{D}q_{1,D}(x)\abs{v_{f}}^{\alpha_{1}}v_{f} - R(x,\abs{v_{f}})v_{f} \text{ in $\Omega$} ,\quad \left.z_{f}\right|_{\partial\Omega}=0. 
\end{equation*}
Multiplying the equation above by $\overline{v_{f}}$ and integrating over $\Omega$ yields 
\begin{equation*}
\int_{\Omega} ((\Delta+q_{0})z_{f}) \overline{v_{f}} \,\rmd x = - \int_{D}q_{1,D}(x)\abs{v_{f}}^{2+\alpha_{1}}\,\rmd x - \int_{\Omega} R(x,\abs{v_{f}})\abs{v_{f}}^{2}\,\rmd x. 
\end{equation*}
The lemma then follows by applying integration by parts twice.
\end{proof}

\subsection{Test data: Complex geometrical optics (CGO) solution for linearized problem}
\addtocontents{toc}{\SkipTocEntry}

The representation formula in \Cref{lem:simplification} and its remark bring us back to the basic mechanism of the classical enclosure method.
Apart from the higher-order remainder term, the principal contribution to the auxiliary indicator functional is given by the first integral over $D$ involving the unknown contrast $q_{1,D}$ and the solution $v_f$ of the linearized problem.
In the classical enclosure method, one chooses suitable special solutions of the background equation and uses their boundary traces as test data to determine whether a prescribed half-space intersects the unknown inclusion $D$.
In the present setting, however, the representation above shows that the relevant solution in the principal term is $v_f$, rather than a solution of the nonlinear background equation.

Since $v_f$ is independent of the unknown inclusion $D$, we may therefore choose $v_f$ to be a suitable special solution of the linearized problem \eqref{eq:linear-part} and use its trace as the test data.
Following Ikehata's enclosure method \cite{Ikehata99Enclosure,Ikehata00Enclosure}, we choose a complex geometrical optics (CGO) solution of the linearized problem, with its boundary trace serving as the test data. Such CGO solutions were first introduced in \cite{sylvester1987global}. We use the standard CGO construction and remainder estimate; see, for example,
\cite[Lemma~2.1]{INUW14increasingstability} for $n\ge3$. The corresponding construction and $H^s$-estimate for $n\ge2$ are also stated in \cite[Remark~2.2]{LLLS21semilinear}. 
For simplicity, we work with a \(C^\infty\) boundary when applying the above CGO construction.
We record below the form needed in our argument. 

\begin{lemma}\label{lem:CGO}
Let $\Omega$ be a bounded $C^\infty$-domain, let $\omega,\omega^\perp\in\mS^{n-1}$ satisfy
$\omega\cdot\omega^\perp=0$, and suppose that
$q_0\in H^s(\Omega)$ for some $s>n/2$ with $s>2$.
Then there exists $h_0>0$ such that, for every $h\in(0,h_0)$,
there exists a solution $v$ to
\[
(\Delta+q_0)v=0\quad\text{in }\Omega
\]
of the form
\[
v=e^{-\frac{\rho\cdot x}{h}}(1+r_\rho),
\qquad
\rho=\omega+\bfi\omega^\perp,
\]
such that
\[
\norm{r_\rho}_{H^s(\Omega)}
\le C h\norm{q_0}_{H^s(\Omega)}
\]
for some constant $C>0$ independent of $h$.
In particular, by Sobolev embedding,
\begin{align}\label{eq:remainder-estimate-A2}
\norm{r_\rho}_{\mA^2(\Omega)}
\le C h\norm{q_0}_{H^s(\Omega)}.    
\end{align}
\end{lemma}

To use the CGO solutions from \Cref{lem:CGO} as test solutions, we need to ensure that their Dirichlet traces belong to $U_{\delta_0}$. We therefore rescale them as follows.

\begin{lemma}\label{lem:CGO-adjusted} 
Let $\Omega$ be a bounded $C^{\infty}$-domain, let $\omega,\omega^{\perp}\in\mS^{n-1}$ with $\omega\cdot\omega^{\perp}=0$ and let $q_{0}\in H^s(\Omega)$ for some $s>n/2$ with $s>2$. Given a parameter $t\in\mR$ and a large parameter $J>0$, we consider the CGO solution 
\begin{equation*}
v_{h} = e^{-\frac{J-t}{h}} e^{-\frac{\rho\cdot x}{h}}(1+r_{\rho}) = e^{-\frac{J}{h}} e^{-\frac{1}{h}(\omega\cdot x - t + \bfi\omega^{\perp}\cdot x)}(1+r_{\rho}) 
\end{equation*}
and we denote its Dirichlet trace $f_{h}:=v_{h}|_{\partial\Omega}$. If we write 
\begin{equation*}
b(\omega) := \inf_{x\in\Omega} \omega\cdot x, 
\end{equation*}
then 
\begin{equation*}
\norm{f_{h}}_{\mA^{\frac{3}{2}}(\partial\Omega)} \le C h^{-2}e^{-\frac{J}{h}}e^{\frac{t-b(\omega)}{h}} \quad \text{for all sufficient small $h\in(0,1)$}
\end{equation*}
for some positive constant $C$ which is independent of the parameter $h$. 
\end{lemma}

\begin{proof} 
For each $j,k=1,\cdots,n$, one computes 
\begin{equation*}
\partial_{j}\partial_{k}v_{h} = e^{-\frac{J}{h}}e^{-\frac{1}{h}(\rho\cdot x-t)}\left(\left(\frac{\rho_{j}}{h}\right)\left(\frac{\rho_{k}}{h}\right)(1+r_{\rho}) - \left(\frac{\rho_{k}}{h}\right)\partial_{j}r_{\rho} - \left(\frac{\rho_{j}}{h}\right)\partial_{k}r_{\rho} + \partial_{j}\partial_{k}r_{\rho} \right), 
\end{equation*}
and then using \eqref{eq:remainder-estimate-A2} we see that 
\begin{equation*}
\norm{\partial_{j}\partial_{k}v_{h}}_{L^{2}(\Omega)} \le C e^{-\frac{J}{h}}e^{\frac{t-b(\omega)}{h}}h^{-2}
\end{equation*}
for some positive constant $C$ which is independent of $h$. Since 
\begin{equation*}
\norm{v_{h}}_{L^{2}(\Omega)} \le C e^{-\frac{J}{h}}e^{\frac{t-b(\omega)}{h}}, 
\end{equation*}
then the trace theorem \cite[Theorem~7.3.11]{AH09TheoreticalNumericalAnalysis} implies 
\begin{equation*}
\norm{f_{h}}_{H^{\frac{3}{2}}(\partial\Omega)}\le C \norm{v_{h}}_{H^{2}(\Omega)} \le C e^{-\frac{J}{h}}e^{\frac{t-b(\omega)}{h}}h^{-2} 
\end{equation*}
for some positive constant $C$ which is independent of $h$. On the other hand, one sees that 
\begin{equation*}
\norm{f_{h}}_{L^{\infty}(\partial\Omega)} \le 2 \sup_{x\in\pl\Omega} \abs{e^{-\frac{J}{h}}e^{-\frac{1}{h}(\omega\cdot x-t+\bfi\omega^{\perp}\cdot x)}} \le 2e^{-\frac{J}{h}}e^{\frac{t-b(\omega)}{h}}.
\end{equation*}
Combining the above two estimates yields
\[
\norm{f_h}_{\mA^{\frac32}(\partial\Omega)}
\le C h^{-2}e^{-\frac{J}{h}}e^{\frac{t-b(\omega)}{h}},
\]
and thus the lemma follows. 
\end{proof}

As an immediate consequence, if $J>B(\omega)-b(\omega)$ and $t<B(\omega)$, where 
\begin{equation*}
B(\omega) := \sup_{x\in\Omega} \omega\cdot x, 
\end{equation*}
then for all sufficiently small $h\in(0,1)$, we have $f_{h}\in U_{\delta_{0}}$. This ensures the CGO solution \Cref{lem:CGO-adjusted} can be substituted into the indicator functionals $\mI$ and $\widetilde{\mI}$. 

\subsection{Reconstruction of the convex hull} 
\addtocontents{toc}{\SkipTocEntry}

We can now finally prove our main theorem as follows:

\begin{theorem}[see also \Cref{thm:main} above]\label{thm:main-detailed}
Let $\Omega \subset \mathbb{R}^n$ ($n \ge 2$) be a bounded $C^{\infty}$-domain and $D$ be an open set satisfying $\overline{D}\subset\Omega$. Assume there exist $0 \ge q_{0}\in H^s(\Omega)$ for some $s>n/2$ with $s>2$, $q_{1}\in L^{\infty}(\Omega)$ with \eqref{eq:perturbation-condition}, a constant $C_{*}>0$, and real numbers $\alpha_{2}>\alpha_{1}>0$, such that the nonlinear term $q$ satisfies \eqref{eq:generalized-Kerr}. Let $t_{*}:=\inf_{x\in D} \omega\cdot x$ and let $f_{h}$ be the function given in \Cref{lem:CGO-adjusted}. For each $\omega\in\mS^{n-1}$, we choose $\omega^\perp\in\mS^{n-1}$ such that $\omega\cdot\omega^\perp=0$ and 
\begin{equation}
J > \max\left\{\frac{(B(\omega)-b(\omega))(2+2\alpha_{1})}{\alpha_{1}} , \frac{(B(\omega)-b(\omega))(2+\alpha_{2})}{\alpha_{2}-\alpha_{1}} \right\} ,\quad b(\omega) < t < B(\omega). \label{eq:choice-parameter}
\end{equation}
Then the following statements hold:
\begin{enumerate}
\renewcommand{\labelenumi}{\theenumi}
\renewcommand{\theenumi}{\rm (\roman{enumi})}
\item \label{itm:1} If $\{\omega\cdot x\le t\}\cap\overline{D}=\emptyset$, that is, $t<t_{*}$, then 
\begin{equation*}
\lim_{h\searrow 0} e^{\frac{(\alpha_{1}+2)J}{h}}\abs{\mI(f_{h})} = 0. 
\end{equation*}
\item \label{itm:2} If $\{\omega\cdot x\le t\}\cap D \neq\emptyset$, that is, $t > t_{*}$, then 
\begin{equation*}
\lim_{h\searrow 0} e^{\frac{(\alpha_{1}+2)J}{h}}\abs{\mI(f_{h})} = +\infty. 
\end{equation*}
\end{enumerate} 
\end{theorem}

\begin{remark*}
\Cref{thm:main-detailed} implies that for any fixed direction $\omega \in \mathbb{S}^{n-1}$, one can determine whether the hyperplane ${ \omega \cdot x = t }$ intersects the inclusion $D$ by testing the limiting behavior of $e^{\frac{(\alpha_{1}+2)J}{h}}\abs{\mathcal{I}(f_h)}$ as $h \to 0$. By varying $t$, we can identify the critical value $t_*$ such that the hyperplane ${ \omega \cdot x = t_* }$ just touches $D$, although the limiting behavior at $t = t_*$ itself is not known. By applying this procedure in every direction $\omega \in \mS^{n-1}$, we reconstruct the convex hull of $D$ and thereby conclude \Cref{thm:main}.  
\end{remark*}

Let us first make some observations before proving \Cref{thm:main-detailed}. Using \Cref{lem:simplification}, we insert the CGO solution $v_{h}$ given in \Cref{lem:CGO-adjusted} into $\widetilde{\mI}$ to see that 
\begin{equation*}
\widetilde{\mI}(f_{h}) = -\int_{D} q_{1,D}(x)\abs{v_{h}}^{2+\alpha_{1}}\,\rmd x - \int_{\Omega} R(x,\abs{v_{h}})\abs{v_{h}}^{2}\,\rmd x. 
\end{equation*}
It follows that 
\begin{equation*}
\begin{aligned} 
e^{\frac{(\alpha_{1}+2)J}{h}}\mI(f_{h}) & = - e^{\frac{(\alpha_{1}+2)J}{h}} \int_{D} q_{1,D}(x)\abs{v_{h}}^{2+\alpha_{1}}\,\rmd x - e^{\frac{(\alpha_{1}+2)J}{h}} \int_{\Omega} R(x,\abs{v_{h}})\abs{v_{h}}^{2}\,\rmd x  \\ 
& \quad + e^{\frac{(\alpha_{1}+2)J}{h}} (\mI(f_{h})-\widetilde{\mI}(f_{h})). 
\end{aligned} 
\end{equation*}
By using \Cref{lem:linearization-indicator}, \Cref{lem:CGO-adjusted} and \eqref{eq:choice-parameter}, we see that 
\begin{equation*}
\begin{aligned}
& e^{\frac{(\alpha_{1}+2)J}{h}}\abs{\mI(f_{h})-\widetilde{\mI}(f_{h})} \le C e^{\frac{(\alpha_{1}+2)J}{h}}\norm{f_{h}}_{\mA^{\frac{3}{2}}(\partial\Omega)}^{2+2\alpha_{1}} \\ 
& \quad \le C h^{-4-4\alpha_{1}} e^{\frac{(\alpha_{1}+2)J}{h}} e^{-\frac{J(2+2\alpha_{1})}{h}} e^{\frac{(t-b(\omega))(2+2\alpha_{1})}{h}} \\ 
& \quad = Ch^{-4-4\alpha_{1}} \exp\left(-\frac{1}{h}(\alpha_{1}J-(t-b(\omega))(2+2\alpha_{1}))\right) \\ 
& \quad \le Ch^{-4-4\alpha_{1}} \exp\left(-\frac{1}{h}(\alpha_{1}J-(B(\omega)-b(\omega))(2+2\alpha_{1}))\right) \\ 
& \quad \rightarrow 0 \quad \text{as $h\searrow 0$.}
\end{aligned}
\end{equation*}
By using \eqref{eq:remainder-estimate-2}, \eqref{eq:remainder-estimate-A2} and \eqref{eq:choice-parameter}, we see that 
\begin{equation*}
\begin{aligned}
& e^{\frac{(\alpha_{1}+2)J}{h}}\left|\int_{\Omega} R(x,\abs{v_{h}})\abs{v_{h}}^{2}\,\rmd x \right| \le C_{*}e^{\frac{(\alpha_{1}+2)J}{h}}\int_{\Omega} \abs{v_{h}}^{2+\alpha_{2}}\,\rmd x \\ 
& \quad \le C_{*} e^{\frac{(\alpha_{1}+2)J}{h}}(1 + \norm{r_{\rho}}_{L^{\infty}(\Omega)})^{2+\alpha_{2}} \int_{\Omega} e^{-\frac{J}{h}(2+\alpha_{2})} e^{-\frac{2+\alpha_{2}}{h}(\omega\cdot x-t)}\,\rmd x \\ 
& \quad \le C \left(1 + C\norm{q_{0}}_{H^s(\Omega)}h \right)^{2+\alpha_{2}} e^{-\frac{J}{h}(\alpha_{2}-\alpha_{1})} e^{-\frac{2+\alpha_{2}}{h}(b(\omega)-B(\omega))} \\ 
& \quad \rightarrow 0 \quad \text{as $h\searrow 0$.}
\end{aligned}
\end{equation*}
Therefore we see that  
\begin{equation}
\lim_{h\searrow 0} e^{\frac{(\alpha_{1}+2)J}{h}}\abs{\mI(f_{h})} = \lim_{h\searrow 0} e^{\frac{(\alpha_{1}+2)J}{h}} \left| \int_{D} q_{1,D}(x)\abs{v_{h}}^{2+\alpha_{1}}\,\rmd x \right|. \label{eq:reduction} 
\end{equation}
We are now ready to prove our main result \Cref{thm:main-detailed}. 

\begin{proof}[Proof of \Cref{thm:main-detailed}\ref{itm:1}] 
Using \Cref{lem:CGO-adjusted} and \eqref{eq:remainder-estimate-A2} we estimate 
\begin{equation*}
\begin{aligned}
& e^{\frac{(\alpha_{1}+2)J}{h}} \left| \int_{D} q_{1,D}(x)\abs{v_{h}}^{2+\alpha_{1}}\,\rmd x \right| \\ 
& \quad \le e^{\frac{(\alpha_{1}+2)J}{h}} \norm{q_{1,D}}_{L^{\infty}(D)} (1+\norm{r_{\rho}}_{L^{\infty}(\Omega)})^{2+\alpha_{1}} \left| \int_{D} e^{-\frac{J}{h}(2+\alpha_{1})}e^{-\frac{2+\alpha_{1}}{h}(\omega\cdot x-t)}\,\rmd x \right| \\ 
& \quad \le C \norm{q_{1,D}}_{L^{\infty}(D)} \left(1 + C\norm{q_{0}}_{H^s(\Omega)}h \right)^{2+\alpha_{1}} e^{-\frac{2+\alpha_{1}}{h}(t_{*}-t)} \\ 
& \quad \rightarrow 0 \quad \text{as $h\searrow 0$,}
\end{aligned}
\end{equation*}
since $t<t_{*}$, hence, in conjunction with \eqref{eq:reduction}, proving the result. 
\end{proof}

\begin{proof}[Proof of \Cref{thm:main-detailed}\ref{itm:2}] 
Using \eqref{eq:contrast-unknown}, we obtain 
\begin{equation*}
\begin{aligned}
& e^{\frac{(\alpha_{1}+2)J}{h}} \left| \int_{D} q_{1,D}(x)\abs{v_{h}}^{2+\alpha_{1}}\,\rmd x \right| \ge e^{\frac{(\alpha_{1}+2)J}{h}} \mu \int_{D} \abs{v_{h}}^{2+\alpha_{1}}\,\rmd x \\ 
& \quad \ge \mu \int_{D} e^{-\frac{2+\alpha_{1}}{h}(\omega\cdot x-t)}\abs{1+r_{\rho}}^{2+\alpha_{1}}\,\rmd x. 
\end{aligned}
\end{equation*}
Note that \eqref{eq:remainder-estimate-A2} ensures that 
\begin{equation*}
\norm{r_{\rho}}_{L^{\infty}(\Omega)} \le  C\norm{q_{0}}_{H^s(\Omega)}h \le \frac{1}{2} \quad \text{for all sufficiently small $h\in(0,1)$,}  
\end{equation*}
therefore we reach 
\begin{equation*}
e^{\frac{(\alpha_{1}+2)J}{h}} \left| \int_{D} q_{1,D}(x)\abs{v_{h}}^{2+\alpha_{1}}\,\rmd x \right| \ge 2^{-2-\alpha_{1}} \mu \int_{D} e^{-\frac{2+\alpha_{1}}{h}(\omega\cdot x-t)} \,\rmd x. 
\end{equation*}
Since $D$ is open, then there exists $\epsilon>0$ such that $t_{*}<t-\epsilon$ and the set $\{\omega\cdot x\le t-\epsilon\}\cap D$ has positive measure. We now have 
\begin{equation*}
\begin{aligned}
& e^{\frac{(\alpha_{1}+2)J}{h}} \left| \int_{D} q_{1,D}(x)\abs{v_{h}}^{2+\alpha_{1}}\,\rmd x \right| \ge 2^{-2-\alpha_{1}} \mu \int_{\{\omega\cdot x\le t-\epsilon\}\cap D} e^{-\frac{2+\alpha_{1}}{h}(\omega\cdot x-t)} \,\rmd x \\ 
& \quad \ge 2^{-2-\alpha_{1}} \mu \int_{\{\omega\cdot x\le t-\epsilon\}\cap D} e^{\frac{2+\alpha_{1}}{h}\epsilon} \,\rmd x = 2^{-2-\alpha_{1}} \mu \abs{\{\omega\cdot x\le t-\epsilon\}\cap D} e^{\frac{2+\alpha_{1}}{h}\epsilon} \\ 
& \quad \rightarrow +\infty \quad \text{as $h\searrow 0$,}
\end{aligned}
\end{equation*}
hence, in conjunction with \eqref{eq:reduction}, proving the result. 
\end{proof}

\section{Numerical implementation}\label{sec:numerical}

A natural way to illustrate \Cref{thm:main-detailed} is to reconstruct an inclusion from simulated measurements, that is, from Neumann data $\Lambda_{q}(f)=\partial_{\nu}u_{f}|_{\partial\Omega}$ produced by a forward solver, in the same way as one would proceed with real data. In this section, we describe and carry out such an implementation. For numerical studies of the enclosure method in the context of electrical impedance tomography, we refer the reader to \cite{IS00EnclosureNumerical,IS04EnclosureNumerical,SRHIS25EnclosureNumerical}. 

The numerical simulations are challenging to carry out using double-precision arithmetic. This difficulty arises from the exponential decay of the CGO test data, which is inherent to the enclosure method, rather than from the nonlinearity of the problem. We carefully address this numerical issue in the design and implementation of \Cref{alg:implementation}. 

\subsection{Setting} 
\addtocontents{toc}{\SkipTocEntry} 

We take $\Omega=B_{1}\subset\mR^{2}$, $q_{0}\equiv0$ and $q(x,\abs{z})=q_{1}(x)\abs{z}^{2}$ with $q_{1,b}\equiv0$ and $q_{1,D}\equiv1$, so that $q_{1}=\chi_{D}$ and the problem \eqref{eq:generalized} reads 
\begin{equation}
\Delta u_{f} + \chi_{D}\abs{u_{f}}^{2}u_{f}=0 \text{ in $B_{1}\subset\mR^{2}$} ,\quad u_{f}|_{\partial B_{1}}=f. \label{eq:nonlinear-PDE}
\end{equation} 
This nonlinearity satisfies \eqref{eq:generalized-Kerr} with $\alpha_{1}=2$ and any $\alpha_{2}>2$. Since $q_{0}\equiv0$ and $q_{1,b}\equiv0$, the background problem is the Laplace equation. Hence the solution $v_{f}$ of \eqref{eq:linear-part}, the background solution $u_{b,f}$ and the background approximate solution $\widetilde{u}_{b,f}$ of \eqref{eq:background-linearized} all coincide: 
\begin{equation*}
v_{f} = u_{b,f} = \widetilde{u}_{b,f} ,\quad \Delta v_{f}=0\text{ in $B_{1}$} ,\quad v_{f}|_{\partial B_{1}}=f. 
\end{equation*}
Consequently, the indicator functional takes the form 
\begin{equation}
\mI(f) = \int_{\partial B_{1}} (\partial_{\nu}u_{f} - \partial_{\nu}v_{f})\overline{f}\,\rmd S. \label{eq:indicator-numerical}
\end{equation}

\subsection{Simulation of the Neumann data}  
\addtocontents{toc}{\SkipTocEntry} 

To simulate the measurement $\Lambda_{q}(f)$, we have to solve the nonlinear problem \eqref{eq:nonlinear-PDE}. 
Let $f\in C^{\infty}(\partial B_{1})$ be complex-valued with $\varepsilon:=\norm{f}_{L^{\infty}(\partial B_{1})}>0$. Since the data of interest are very small, we separate their size from their shape. We set $\hat{f}=\varepsilon^{-1}f$ and let $v_{\hat{f}}$ be the harmonic extension of $\hat{f}$, so that $v_{f}=\varepsilon v_{\hat{f}}$ and, by the maximum principle, $\norm{v_{\hat{f}}}_{L^{\infty}(B_{1})}\le1$. The nonlinear part of the solution is captured by the rescaled correction 
\begin{equation*}
W := \varepsilon^{-3}\left( u_{f}-\varepsilon v_{\hat{f}} \right).  
\end{equation*}
A direct computation shows that $W$ satisfies 
\begin{subequations} \label{eq:Picard}
\begin{equation}
\Delta W = -\chi_{D}g(\varepsilon^{2}W + v_{\hat{f}}) \text{ in $B_{1}$} ,\quad W|_{\partial B_{1}}=0, \label{eq:Picard-limit} 
\end{equation}
where $g(z):=\abs{z}^{2}z$. 
Motivated by the fixed-point argument in the proof of \Cref{lem:well-posedness}, we approximate $W$ by the Picard iteration 
\begin{equation}
\Delta W^{(\ell)} = -\chi_{D}g(\varepsilon^{2} W^{(\ell-1)} + v_{\hat{f}}) \text{ in $B_{1}$} ,\quad W^{(\ell)}|_{\partial B_{1}}=0,\quad \ell\in\mN, \label{eq:Picard-iteration}  
\end{equation} 
\end{subequations} 
with $W^{(0)}=0$. After $\ell^{*}\in\mN$ steps, the simulated measurement is 
\begin{equation}
\Lambda^{(\ell^{*})}(f) := \left.\partial_{\nu}\left( \varepsilon v_{\hat{f}} + \varepsilon^{3}W^{(\ell^{*})} \right)\right|_{\partial B_{1}}. \label{eq:simulated-measurement} 
\end{equation}
Note that the first iterate reproduces the approximate solution of \eqref{eq:approximate-solution}, namely $\tilde{u}_{f}=\varepsilon v_{\hat{f}} + \varepsilon^{3}W^{(1)}$. 

The following lemma shows that the iteration is well defined and converges under a mild condition on $\varepsilon$, which, unlike \Cref{lem:well-posedness}, does not depend on any Sobolev norm of $f$. 

\begin{lemma}\label{lem:well-posed-relax}
Assume $\varepsilon^{2}\le1/2$, and let 
\begin{equation*}
\mathscr{X}:=\left\{W\in L^{\infty}(B_{1}) : \norm{W}_{L^{\infty}(B_{1})}\le1\right\}. 
\end{equation*}
Then \eqref{eq:Picard-limit} has a unique solution $W\in\mathscr{X}$, all iterates $W^{(\ell)}$ belong to $\mathscr{X}$, and 
\begin{equation*}
\norm{W-W^{(\ell)}}_{L^{\infty}(B_{1})} \le \kappa^{\ell} \quad \text{for all $\ell\in\mN$,} 
\end{equation*}
where $\kappa:=\frac{27}{16}\varepsilon^{2} \le \frac{27}{32}<1$. In particular, \eqref{eq:nonlinear-PDE} has a unique solution $u_{f}$ satisfying 
\begin{equation*}
\norm{u_{f}-v_{f}}_{L^{\infty}(B_{1})}\le\varepsilon^{3}. 
\end{equation*}
\end{lemma}

\begin{proof}
Let $\mathcal{S}_{0}$ denote the solution operator of $\Delta w=F$ in $B_{1}$ with $w|_{\partial B_{1}}=0$. Let $\vartheta\in\mR$ be a number to be chosen later. Comparing $\Re(e^{\bfi\vartheta}\mathcal{S}_{0}F)$ with $\norm{F}_{L^{\infty}(B_{1})}(1-\abs{x}^{2})/4$, the maximum principle gives 
\begin{equation*}
\abs{\Re(e^{\bfi\vartheta}(\mathcal{S}_{0}F)(x))} \le \frac{1}{4}\norm{F}_{L^{\infty}(B_{1})} \quad \text{for all $\vartheta\in\mR$ and a.e. $x\in B_{1}$.}
\end{equation*}
For a.e. $x\in B_{1}$, choosing $\vartheta=-\arg(\mathcal{S}_{0}F)$ yields 
\begin{equation*}
\abs{(\mathcal{S}_{0}F)(x)}\le\frac{1}{4}\norm{F}_{L^{\infty}(B_{1})} \quad \text{for a.e. $x\in B_{1}$,} 
\end{equation*}
that is, 
\begin{equation}
\norm{\mathcal{S}_{0}F}_{L^{\infty}(B_{1})} \le \frac{1}{4}\norm{F}_{L^{\infty}(B_{1})}. \label{eq:bound-S0}
\end{equation}
Let $W\in \mathscr{X}$. Since $\norm{v_{\hat{f}}}_{L^{\infty}(B_{1})}\le1$, $\norm{W}_{L^{\infty}(B_{1})}\le1$ and $\varepsilon^{2}\le\frac{1}{2}$, we have pointwise 
\begin{equation}
\norm{\varepsilon^{2}W + v_{\hat{f}}}_{L^{\infty}(B_{1})} \le \varepsilon^{2}\norm{W}_{L^{\infty}(B_{1})} + \norm{v_{\hat{f}}}_{L^{\infty}(B_{1})} \le \frac{1}{2}+1= \frac{3}{2}. \label{eq:contraction1-extend}
\end{equation}
Hence $g(\varepsilon^{2}W+v_{\hat{f}})$ satisfies 
\begin{equation*}
\norm{g(\varepsilon^{2}W+v_{\hat{f}})}_{L^{\infty}(B_{1})} \le \norm{\varepsilon^{2}W+v_{\hat{f}}}_{L^{\infty}(B_{1})}^{3}\le \left(\frac{3}{2}\right)^{3}=\frac{27}{8}. 
\end{equation*}
Choosing $F=\chi_{D}g(\varepsilon^{2}W+v_{\hat{f}})$ in \eqref{eq:bound-S0}, we have 
\begin{equation*}
\norm{\mathcal{S}_{0}(\chi_{D}g(\varepsilon^{2}W+v_{\hat{f}}))}_{L^{\infty}(B_{1})} \le \frac{1}{4}\cdot\frac{27}{8} = \frac{27}{32} < 1, 
\end{equation*}
which shows that $\mathcal{S}_{0}(\chi_{D}g(\varepsilon^{2}W+v_{\hat{f}}))\in\mathscr{X}$. 

We see that the real derivative of $g$ is $Dg(z)\zeta=\abs{z}^{2}\zeta+2\Re(\overline{z}\zeta)z$ for all $\zeta\in\mC$, so that 
\begin{equation}
\abs{Dg(z)\zeta} \le 3\abs{z}^{2}\abs{\zeta} \quad \text{for all $\zeta\in\mC$.} \label{eq:gradient-g}
\end{equation}
Let $W_{1},W_{2}\in\mathscr{X}$ . By \eqref{eq:contraction1-extend}, we see that 
\begin{equation*}
\norm{\varepsilon^{2}W_{j}+v_{\hat{f}}}_{L^{\infty}(B_{1})} \le \frac{3}{2} \quad \text{for all $j=1,2$.} 
\end{equation*}
The disk $\{z\in\mC:\abs{z}\le\frac{3}{2}\}$ is convex, and the mean value inequality gives 
\begin{equation}
\begin{aligned}
& \norm{g(\varepsilon^{2}W_{1}+v_{\hat{f}})-g(\varepsilon^{2}W_{2}+v_{\hat{f}})}_{L^{\infty}(B_{1})} \\
&\quad \le 3\left(\frac{3}{2}\right)^{2}\varepsilon^{2}\norm{W_{1}-W_{2}}_{L^{\infty}(B_{1})} = \frac{27}{4}\varepsilon^{2}\norm{W_{1}-W_{2}}_{L^{\infty}(B_{1})}. 
\end{aligned} \label{eq:MVT-consequence}
\end{equation}
Choosing $F=\chi_{D}g(\varepsilon^{2}W_{1}+v_{\hat{f}})-\chi_{D}g(\varepsilon^{2}W_{2}+v_{\hat{f}})$ in \eqref{eq:bound-S0}, we have 
\begin{equation*}
\begin{aligned}
& \norm{\mathcal{S}_{0}(\chi_{D}g(\varepsilon^{2}W_{1}+v_{\hat{f}})) - \mathcal{S}_{0}(\chi_{D}g(\varepsilon^{2}W_{2}+v_{\hat{f}}))}_{L^{\infty}(B_{1})} \\ 
&\quad \le \frac{1}{4}\cdot\frac{27}{4}\varepsilon^{2}\norm{W_{1}-W_{2}}_{L^{\infty}(B_{1})} = \kappa \norm{W_{1}-W_{2}}_{L^{\infty}(B_{1})}. 
\end{aligned}
\end{equation*}
Therefore, the map 
\begin{equation}
W\mapsto -\mathcal{S}_{0}\left(\chi_{D}g(\varepsilon^{2}W + v_{\hat{f}})\right) \label{eq:fix-map}
\end{equation}
sends $\mathscr{X}$ into itself and is a contraction on $\mathscr{X}$. The Banach fixed-point theorem then guarantees that \eqref{eq:Picard-limit} has a unique solution $W\in\mathscr{X}$. 

Since $W^{(0)}\equiv0$, the above argument shows that every $W^{(\ell)}$ belongs to $\mathscr{X}$ and satisfies 
\begin{equation*}
\norm{W-W^{(\ell)}}_{L^{\infty}(B_{1})} \le \kappa^{\ell} \quad \text{for all $\ell\in\mN$.}  
\end{equation*}
Finally, note that $u_{f}:=\varepsilon v_{\hat{f}}+\varepsilon^{3}W$ solves \eqref{eq:nonlinear-PDE} and satisfies 
\begin{equation*}
\norm{u_{f}-v_{f}}_{L^{\infty}(B_{1})} = \varepsilon^{3}\norm{W}_{L^{\infty}(B_{1})} \le \varepsilon^{3}. 
\end{equation*}
Conversely, if $u$ solves \eqref{eq:nonlinear-PDE} with $\norm{u-v_{f}}_{L^{\infty}(B_{1})}\le\varepsilon^{3}$, then $\varepsilon^{-3}(u-v_{f})\in\mathscr{X}$ is a fixed point on the map \eqref{eq:fix-map}, which is equal to $W$. This proves uniqueness. 
\end{proof}

\subsection{Reconstruction from Neumann data}   
\addtocontents{toc}{\SkipTocEntry} 

Inserting the simulated measurement \eqref{eq:simulated-measurement} into \eqref{eq:indicator-numerical}, we obtain the simulated indicator 
\begin{equation*}
\mI^{(\ell^{*})}(f) := \int_{\partial B_{1}} \left( \Lambda^{(\ell^{*})}(f) - \partial_{\nu}v_{f} \right) \overline{f}\,\rmd S. 
\end{equation*}
For the analysis below, it is useful to rewrite the indicator in a form that does not involve a difference of Neumann data. Since $W|_{\partial B_{1}}=0$ and $\Delta\overline{v_{\hat{f}}}=0$, integration by parts and \eqref{eq:Picard-limit} gives 
\begin{subequations} \label{eq:indicator-numerical-analysis}
\begin{equation}
\mI(f) = \varepsilon^{4}\int_{\partial B_{1}}\partial_{\nu}W \overline{\hat{f}}\,\rmd S = \varepsilon^{4}\int_{B_{1}} \Delta W \overline{v_{\hat{f}}} \,\rmd x = -\varepsilon^{4}\int_{D} g\left(\varepsilon^{2}W + v_{\hat{f}}\right)\overline{v_{\hat{f}}}\,\rmd x, \label{eq:indicator-numerical-1}
\end{equation}
and, in the same way, 
\begin{equation}
\mI^{(\ell^{*})}(f) = -\varepsilon^{4} \int_{D} g\left(\varepsilon^{2}W^{(\ell^{*}-1)} + v_{\hat{f}}\right)\overline{v_{\hat{f}}}\,\rmd x. \label{eq:indicator-numerical-analysis-approx}
\end{equation} 
\end{subequations} 

\begin{remark*} 
For $\ell^{*}=1$, this is  
\begin{equation}
\mI^{(1)}(f) = -\varepsilon^{4}\int_{D}\abs{v_{\hat{f}}}^{4}\,\rmd x = \widetilde{\mI}(f), \label{eq:I1-special-repn}
\end{equation}
in accordance with \Cref{lem:simplification}, since $R(x,\abs{z})\equiv0$. 
\end{remark*} 

We stress that \eqref{eq:indicator-numerical-analysis} involve the inclusion $D$ explicitly. They are tools for analyzing the size of the indicator, not substitutes for the measurement. Evaluating the indicator through them would amount to an ``inverse crime''. 

The error caused by stopping the iteration is controlled as follows. Subtracting \eqref{eq:indicator-numerical-analysis-approx} from \eqref{eq:indicator-numerical-1}, and using $\norm{v_{\hat{f}}}_{L^{\infty}(B_{1})}\le1$, \eqref{eq:MVT-consequence} and \Cref{lem:well-posed-relax}, we obtain 
\begin{equation*}
\begin{aligned}
& \abs{\mI(f)-\mI^{(\ell^{*})}(f)} \le \varepsilon^{4}\int_{D}\abs{g(\varepsilon^{2}W+v_{\hat{f}})-g(\varepsilon^{2}W^{(\ell^{*}-1)}+v_{\hat{f}})}\,\rmd x \\ 
&\quad \le \frac{27}{4}\abs{D}\varepsilon^{6}\kappa^{\ell^{*}-1} = \frac{27}{4}\left(\frac{27}{16}\right)^{\ell^{*}-1}\abs{D}\varepsilon^{4+2\ell^{*}}. 
\end{aligned}
\end{equation*}
Thus, the simulated indicator  converges to the true one as $\ell^{*}\rightarrow\infty$.

\subsection{Reconstruction procedure}    
\addtocontents{toc}{\SkipTocEntry} 

For $\omega=(\cos\theta,\sin\theta)\in\mS^{1}$ we set $\omega^{\perp}:=(-\sin\theta,\cos\theta)$ and $\rho:=\omega+\bfi\omega^{\perp}$, so that $\rho\cdot\rho=0$. Since $q_{0}\equiv0$, the CGO solution in \Cref{lem:CGO} is exact, i.e. $r_{\rho}\equiv0$, and the test data is 
\begin{equation*}
f_{h}(\hat{x}) = e^{-\frac{J}{h}} e^{-\frac{1}{h}(\rho\cdot\hat{x}-t)} \quad \text{for all $\hat{x}\in\partial B_{1}$,} 
\end{equation*}
with $b(\omega)=-1<t<1=B(\omega)$. Since $\abs{f_{h}(\hat{x})}=e^{-\frac{J}{h}}e^{-\frac{1}{h}(\omega\cdot\hat{x}-t)}$ attains its maximum at $\hat{x}=-\omega$, we have $\varepsilon = \norm{f_{h}}_{L^{\infty}(\partial B_{1})} = e^{-\frac{1}{h}(J-t-1)}$. 

If $J$ were fixed as required by \Cref{thm:main-detailed}, say $J=10$ with $h=1/10$, then $\varepsilon\le \exp(-80)$. The nonlinear part of the Neumann data would then be smaller than the background part by a factor of at most $\varepsilon^{2}\le\exp(-160) = O(10^{-70})$, which rules out any computation in floating-point arithmetic from the outset. 

To give the method the best possible chance, we therefore relax this requirement. Observe that $f_{h}$ depends on $(J,t)$ only through $J-t$. Instead of fixing $J$, we fix the size $\varepsilon_{0}\in(0,1/\sqrt{2}]$ of the data, and for each $t\in(-1,1)$ we choose 
\begin{equation}
J = J_{t} := t+1+h\log\frac{1}{\varepsilon_{0}} ,\quad \text{so that $\varepsilon=\varepsilon_{0}$ and $f_{h}(\hat{x})=\varepsilon_{0}e^{-\frac{1}{h}(1+\rho\cdot\hat{x})}$.} \label{eq:choice-J}
\end{equation}
Thus the boundary datum is the same for every $t$, and its normalized harmonic extension is known in closed form: 
\begin{equation*}
v_{\hat{f}_{h}}(x) = e^{-\frac{1}{h}(1+\rho\cdot x)} \quad \text{for all $x\in B_{1}$.} 
\end{equation*}
We consider $h=1/20$. Since the size of the indicator relative to the Neumann data grows with $\varepsilon_{0}$, we take the largest value admitted by \Cref{lem:well-posed-relax}, $\varepsilon_{0}=1/\sqrt{2}$, so that $J_{t}=t+1+\frac{1}{40}\log 2\approx t+1.017$.  

With this choice, the indicator takes a particularly simple form. Writing $W$ and $W^{(\ell)}$ for the solutions of \eqref{eq:Picard} with $\varepsilon=\varepsilon_{0}$ and $v_{\hat{f}}=v_{\hat{f}_{h}}$, we set 
\begin{equation*}
K(\omega) := \int_{D}g\left( \varepsilon_{0}^{2}W+v_{\hat{f}_{h}} \right)\overline{v_{\hat{f}_{h}}}\,\rmd x ,\quad K^{(\ell)}(\omega) := \int_{D}g\left( \varepsilon_{0}^{2}W^{(\ell-1)}+v_{\hat{f}_{h}} \right)\overline{v_{\hat{f}_{h}}}\,\rmd x, 
\end{equation*}
which do not depend on $t$. 

\begin{remark*}
For $\ell^{*}=1$, this is 
\begin{equation*}
K^{(1)}(\omega) = \int_{D}\abs{v_{\hat{f}_{h}}}^{4}\,\rmd x. 
\end{equation*}
\end{remark*}

Since $e^{4J_{t}/h}\varepsilon_{0}^{4}=e^{4(t+1)/h}$, it follows from \eqref{eq:indicator-numerical-1} that the normalized indicator of \Cref{thm:main-detailed} is 
\begin{subequations}\label{eq:indicators}
\begin{equation}
\Phi(\omega,t) := \log\left(e^{\frac{(\alpha_{1}+2)J_{t}}{h}}\abs{\mI(f_{h})}\right) = \frac{4(t+1)}{h} + \log\abs{K(\omega)}. \label{eq:indicators-1}  
\end{equation}
Similarly, the simulated normalized indicator is given by 
\begin{equation}
\Phi^{(\ell^{*})}(\omega,t) := \log\left(e^{\frac{(\alpha_{1}+2)J_{t}}{h}}\abs{\mI^{(\ell^{*})}(f_{h})}\right) = \frac{4(t+1)}{h} + \log\abs{K^{(\ell^{*})}(\omega)}. \label{eq:indicators-2}  
\end{equation}
\end{subequations} 
Both $t\mapsto\Phi(\omega,t)$ and $t\mapsto\Phi^{(\ell^{*})}(\omega,t)$ are affine with slope $4/h$. Their zeros are 
\begin{equation*}
t_{\omega} = -1-\frac{h}{4}\log\abs{K(\omega)} \quad\text{and}\quad t_{\omega}^{(\ell^{*})} = -1-\frac{h}{4}\log\abs{K^{(\ell^{*})}(\omega)}.  
\end{equation*}

However, the choice \eqref{eq:choice-J} does not satisfy condition \eqref{eq:choice-parameter} of \Cref{thm:main-detailed}, which requires $J>6$ independently of $t$. In the proof of \Cref{thm:main-detailed}, that condition serves two purposes, and we now check that the conclusion of the theorem survives without it. 

The first purpose is to ensure that the data are small. By \Cref{lem:CGO-adjusted}, 
\begin{equation*}
\norm{f_{h}}_{\mathcal{A}^{3/2}(\partial B_{1})} \le Ch^{-2}e^{-\frac{1}{h}(J-t+b(\omega))}, 
\end{equation*}
and the condition $J>B(\omega)-b(\omega)$ guarantees $f_{h}\in U_{\delta_{0}}$ for small $h$. With $J=J_{t}$, we only have $\norm{f_{h}}_{\mathcal{A}^{3/2}(\partial B_{1})} \le Ch^{-2}\varepsilon_{0}$, so \Cref{lem:well-posedness} cannot be applied directly. However, \Cref{lem:well-posed-relax} requires only $\varepsilon_{0}^{2}\le\frac{1}{2}$. It therefore provides, for every $h$, a unique solution $u_{f_{h}}$ of \eqref{eq:nonlinear-PDE} close to $v_{f_{h}}$, together with the convergence of the Picard iteration \eqref{eq:Picard}. 

The second purpose is to make the linearization error negligible. In the proof of \Cref{thm:main-detailed}, \Cref{lem:linearization-indicator} gives 
\begin{equation*}
e^{\frac{(\alpha_{1}+2)J}{h}}\abs{\mI(f_{h})-\widetilde{\mI}(f_{h})} \le Ce^{\frac{4J}{h}}\norm{f_{h}}_{\mathcal{A}^{3/2}(\partial B_{1})}^{6}. 
\end{equation*}
With $J=J_{t}$ this general bound only yields 
\begin{equation*}
e^{\frac{4J_{t}}{h}}\abs{\mI(f_{h})-\widetilde{\mI}(f_{h})} \le Ch^{-12}\varepsilon_{0}^{2}e^{\frac{4(t+1)}{h}}, 
\end{equation*}
which is useless. The bound is far from sharp, however, because it ignores the fact that $v_{\hat{f}_{h}}$ is exponentially small in $D$. Exploiting this, we obtain the following sharper estimate. It shows that the linearization error is negligible relative to the principal term, which is all that the proof of \Cref{thm:main-detailed} requires. 

\begin{theorem}\label{thm:approximation}
Let $\varepsilon_{0}\in(0,1/\sqrt{2}]$ and $\omega\in\mS^{1}$, let $t_{*}:=\inf_{x\in D}\omega\cdot x>-1$, and set 
\begin{equation*}
\sfm:=\norm{v_{\hat{f}_{h}}}_{L^{\infty}(D)}=e^{-(1+t_{*})/h}\le1. 
\end{equation*}
Then the following holds: 
\begin{enumerate}
\renewcommand{\labelenumi}{\theenumi}
\renewcommand{\theenumi}{\rm (\alph{enumi})}
\item\label{itm:a} $\norm{W}_{L^{\infty}(B_{1})} \le \sfm^{3}/2$ and $\norm{W^{(\ell)}}_{L^{\infty}(B_{1})}\le\sfm^{3}/2$ for all $\ell\in\mN$. 
\item\label{itm:b} We have 
\begin{equation*}
e^{\frac{4J_{t}}{h}}\abs{\mI(f_{h})-\widetilde{\mI}(f_{h})} = e^{\frac{4(t+1)}{h}}\abs{K(\omega)-K^{(1)}(\omega)} \le 3\varepsilon_{0}^{2}\abs{D}e^{\frac{4(t+1)-6(1+t_{*})}{h}}, 
\end{equation*}
and, for every $\ell^{*}\in\mN$, the same bound holds with $K^{(\ell^{*})}$ in place of $K$. 
\item\label{itm:c} Write $K(\omega) = K^{(1)}(\omega)(1+\delta(\omega))$ with $\delta(\omega)\in\mC$. For every $\eta\in(0,(1+t_{*})/2)$ there exists $c_{\eta}>0$, independent of $h$, such that 
\begin{equation*}
\abs{\delta(\omega)} \le \frac{3\varepsilon_{0}^{2}\abs{D}}{c_{\eta}}e^{-\frac{2(1+t_{*})-4\eta}{h}}. 
\end{equation*}
In particular, 
\begin{equation*}
\left| \frac{\abs{K(\omega)}}{K^{(1)}(\omega)} - 1 \right| \le \abs{\delta(\omega)},  
\end{equation*}
so that $\abs{K(\omega)}/K^{(1)}(\omega)\rightarrow 1$ as $h\searrow0$. 
Likewise, for every $\ell^{*}\in\mN$, if we write $K^{(\ell^{*})}(\omega) = K^{(1)}(\omega)(1+\delta^{(\ell^{*})}(\omega))$ with $\delta^{(\ell^{*})}(\omega)\in\mC$, then $\abs{\delta^{(\ell^{*})}(\omega)}$ satisfies the same bound as $\abs{\delta(\omega)}$, and 
\begin{equation*}
\left| \frac{\abs{K^{(\ell^{*})}(\omega)}}{K^{(1)}(\omega)} - 1 \right| \le \abs{\delta^{(\ell^{*})}(\omega)}. 
\end{equation*} 
\item\label{itm:d} Let $\Phi$ and $\Phi^{(\ell^{*})}$ be given by \eqref{eq:indicators}. As $h\searrow0$, both $\Phi(\omega,t)$ and, for every $\ell^{*}\in\mN$, $\Phi^{(\ell^{*})}(\omega,t)$ are equal to 
\begin{equation*}
\frac{4t}{h} + \log\int_{D} e^{-\frac{4}{h}\omega\cdot x}\,\rmd x + o(1). 
\end{equation*} 
Consequently, both $\Phi(\omega,t)$ and $\Phi^{(\ell^{*})}(\omega,t)$ tend to $-\infty$ if $t<t_{*}$ and to $+\infty$ if $t>t_{*}$. 
\end{enumerate} 
\end{theorem}

\begin{proof}[Proof of \Cref{thm:approximation}\ref{itm:a}] 
By \eqref{eq:bound-S0}, we have 
\begin{equation*}
\begin{aligned} 
&\norm{W}_{L^{\infty}(B_{1})} \le \frac{1}{4}\norm{\chi_{D}g(\varepsilon_{0}^{2}W+v_{\hat{f}_{h}})}_{L^{\infty}(B_{1})} \\
&\quad \le \frac{1}{4}\norm{\varepsilon_{0}^{2}W+v_{\hat{f}_{h}}}_{L^{\infty}(D)}^{3} \le \frac{1}{4}\left( \sfm + \varepsilon_{0}^{2}\norm{W}_{L^{\infty}(B_{1})} \right)^{3} , 
\end{aligned} 
\end{equation*}
that is, the number $s:=\norm{W}_{L^{\infty}(B_{1})}\in[0,1]$ satisfies 
\begin{equation*}
s \le \phi(s) := \frac{1}{4}(\sfm + \varepsilon_{0}^{2}s)^{3}. 
\end{equation*}
The function $\phi(s)-s$ is convex and positive at $s=0$. It is negative at $s=1$, since $\phi(1)\le\frac{1}{4}(\frac{3}{2})^{3}=\frac{27}{32}<1$, and at $s=\sfm^{3}/2$, since 
\begin{equation*}
\phi(\sfm^{3}/2) \le \frac{\sfm^{3}}{4}\left(1+\frac{\sfm^{2}}{4}\right)^{3} < \frac{\sfm^{3}}{2}. 
\end{equation*}
By convexity, $\phi(s)-s<0$ on $[\sfm^{3}/2,1]$, hence $s<\sfm^{3}/2$. For each $\ell\in\mN$, we have 
\begin{equation*}
\norm{W^{(\ell)}}_{L^{\infty}(B_{1})} \le \phi(\norm{W^{(\ell-1)}}_{L^{\infty}(B_{1})}), 
\end{equation*}
so the bound follows by induction on $\ell$, starting from $W^{(0)}\equiv0$. 
\end{proof} 

\begin{proof}[Proof of \Cref{thm:approximation}\ref{itm:b}] 
The identity 
\begin{equation*}
e^{\frac{4J_{t}}{h}}\abs{\mI(f_{h})-\widetilde{\mI}(f_{h})} = e^{\frac{4(t+1)}{h}}\abs{K(\omega)-K^{(1)}(\omega)} 
\end{equation*}
follows from \eqref{eq:indicator-numerical-analysis}, \eqref{eq:I1-special-repn} and $e^{4J_{t}/h}\varepsilon_{0}^{4} = e^{4(t+1)/h}$.

By \Cref{thm:approximation}\ref{itm:a} and $\varepsilon_{0}^{2}\le1/2$, we have 
\begin{equation*} 
\norm{\varepsilon_{0}^{2}W}_{L^{\infty}(B_{1})}\le\sfm^{3}/4. 
\end{equation*} 
Using \eqref{eq:gradient-g} and the mean value inequality gives 
\begin{equation*}
\norm{g(\varepsilon_{0}^{2}W+v_{\hat{f}_{h}}) - g(v_{\hat{f}_{h}})}_{L^{\infty}(D)} \le 3\left(\frac{5}{4}\sfm\right)^{2}\varepsilon_{0}^{2}\frac{\sfm^{3}}{2} \le 3\varepsilon_{0}^{2}\sfm^{5},  
\end{equation*}
because every ``point'' on the ``segment'' joining $\varepsilon_{0}^{2}W+v_{\hat{f}_{h}}$ and $v_{\hat{f}_{h}}$ has the form 
\begin{equation*}
s\varepsilon_{0}^{2}W+v_{\hat{f}_{h}} \quad \text{for some $s\in[0,1]$,} 
\end{equation*}
so that $\norm{s\varepsilon_{0}^{2}W+v_{\hat{f}_{h}}}_{L^{\infty}(D)}\le\norm{\varepsilon_{0}^{2}W}_{L^{\infty}(D)}+\sfm\le\frac{5}{4}\sfm$. Therefore, we obtain 
\begin{equation*}
\abs{K(\omega)-K^{(1)}(\omega)} \le \int_{D} 3\varepsilon_{0}^{2}\sfm^{5}\abs{v}\,\rmd x \le 3\varepsilon_{0}^{2}\abs{D}\sfm^{6} = 3\varepsilon_{0}^{2}\abs{D}e^{-\frac{6(1+t_{*})}{h}}. 
\end{equation*}
Multiplying by $e^{4(t+1)/h}$ gives the claimed bound. For $K^{(\ell^{*})}$, the same argument applies with $W^{(\ell^{*}-1)}$ in place of $W$, since, by \Cref{thm:approximation}\ref{itm:a}, $\norm{W^{(\ell^{*}-1)}}_{L^{\infty}(B_{1})}\le\sfm^{3}/2$ as well. 
\end{proof}

\begin{proof}[Proof of \Cref{thm:approximation}\ref{itm:c}] 
Let $c_{\eta}:=\abs{\left\{x\in\mR^{2}:\omega\cdot x<t_{*}+\eta\right\}\cap D}$, which is positive since $D$ is open. Since $\abs{v(x)}^{4} = e^{-4(1+\omega\cdot x)/h}$, we have 
\begin{equation*} 
K^{(1)}(\omega) \ge c_{\eta}e^{-\frac{4(1+t_{*}+\eta)}{h}} > 0. 
\end{equation*} 
By definition, $\delta(\omega)=(K(\omega)-K^{(1)}(\omega))/K^{(1)}(\omega)$, so \Cref{thm:approximation}\ref{itm:b} gives 
\begin{equation*}
\abs{\delta(\omega)} \le \frac{3\varepsilon_{0}^{2}\abs{D}e^{-\frac{6(1+t_{*})}{h}}}{c_{\eta}e^{-\frac{4(1+t_{*}+\eta)}{h}}} = \frac{3\varepsilon_{0}^{2}\abs{D}}{c_{\eta}}e^{-\frac{2(1+t_{*})-4\eta}{h}}, 
\end{equation*}
which is the claimed bound. 
It remains to compare $\abs{K(\omega)}$ with $K^{(1)}(\omega)$. Since $K^{(1)}(\omega)$ is real and positive, we have 
\begin{equation*}
\frac{\abs{K(\omega)}}{K^{(1)}(\omega)} = \abs{1+\delta(\omega)}. 
\end{equation*}
Then 
\begin{equation*}
\left| \frac{\abs{K(\omega)}}{K^{(1)}(\omega)} - 1 \right| = \left| \abs{1+\delta(\omega)} - 1 \right| \le \abs{\delta(\omega)}. 
\end{equation*}
The same argument, using the bound on $\abs{K^{(\ell^{*})}(\omega)-K^{(1)}(\omega)}$ from \Cref{thm:approximation}\ref{itm:b}, applies to $\delta^{(\ell^{*})}$. 
\end{proof} 

\begin{proof}[Proof of \Cref{thm:approximation}\ref{itm:d}] 
By \Cref{thm:approximation}\ref{itm:c}, $\abs{\delta(\omega)}<1$ for all sufficiently small $h$, in particular $K(\omega)\neq0$, so that $\Phi(\omega,t)$ is well-defined. Then 
\begin{equation*}
\log\abs{K(\omega)} = \log K^{(1)}(\omega) + \log\abs{1+\delta(\omega)}. 
\end{equation*}
Note that 
\begin{equation*}
\abs{\log\abs{1+\delta(\omega)}} \le-\log(1-\abs{\delta(\omega)})\rightarrow0\text{ as $h\searrow0$}
\end{equation*}
because $1-\abs{\delta}\le\abs{1+\delta}\le1+\abs{\delta}$ and $\log(1+\abs{\delta})\le-\log(1-\abs{\delta})$. Since 
\begin{equation*}
K^{(1)}(\omega) = e^{-\frac{4}{h}}\int_{D}e^{-\frac{4}{h}\omega\cdot x}\,\rmd x, 
\end{equation*}
the expansion follows from \eqref{eq:indicators-1}. If $t<t_{*}$, the right-hand side is at most 
\begin{equation*}
\frac{4(t-t_{*})}{h}+\log\abs{D}+o(1)\rightarrow-\infty \quad \text{as $h\searrow0$.} 
\end{equation*}
If $t>t_{*}$, we choose $\eta<t-t_{*}$ in the definition of $c_{\eta}$, and the right-hand side is at least 
\begin{equation*}
\frac{4(t-t_{*}-\eta)}{h} + \log c_{\eta} + o(1) \rightarrow+\infty \quad \text{as $h\searrow0$.} 
\end{equation*}
The same argument, with $\delta^{(\ell^{*})}$ in place of $\delta$ and \eqref{eq:indicators-2} in place of \eqref{eq:indicators-1}, applies to $\Phi^{(\ell^{*})}$. 
\end{proof}

With these choices of test data, the reconstruction procedure is summarized in \Cref{alg:implementation}.

\begin{algorithm}[H] 
\caption{Enclosure method based on simulated Neumann data\label{alg:implementation}} 
\begin{algorithmic}[1] 
\REQUIRE $D$ (for the simulation only), $h$, $\varepsilon_{0},\ell^{*}$, number of directions $N_{\omega}$. 
\FOR{$k=0,\cdots,N_{\omega}-1$}
\STATE $\theta_{k} = 2\pi k/N_{\omega}$, $\omega=(\cos\theta_{k},\sin\theta_{k})$, $\omega^{\perp}=(-\sin\theta_{k},\cos\theta_{k})$, $\rho=\omega+\bfi\omega^{\perp}$ 
\STATE $f(\hat{x})=\varepsilon_{0}e^{-(1+\rho\cdot\hat{x})/h}$ for all $\hat{x}\in\partial B_{1}$ 
\STATE $v_{\hat{f}}(x)=e^{-(1+\rho\cdot x)/h}$ and $v_{f}(x)=\varepsilon_{0}v_{\hat{f}}(x)$ for all $x\in B_{1}$ 
\STATE $W^{(0)}\equiv0$ 
\FOR{$\ell=1,\cdots,\ell^{*}$} 
\STATE Solve \eqref{eq:Picard-iteration} with $\varepsilon=\varepsilon_{0}$ for $W^{(\ell)}$ 
\ENDFOR 
\STATE $\Lambda^{(\ell^{*})}(f) = \partial_{\nu}(\varepsilon_{0}v_{\hat{f}}+\varepsilon_{0}^{3}W^{(\ell^{*})})$ on $\partial B_{1}$ \COMMENT{simulated measurement} 
\STATE $\mI^{(\ell^{*})}(f)=\int_{\partial B_{1}}(\Lambda^{(\ell^{*})}(f) - \partial_{\nu}v_{f})\overline{f}\,\rmd S$ \COMMENT{difference of Neumann data} 
\STATE $t_{\omega_{k}}^{(\ell^{*})}=-1-\frac{h}{4}\log(\varepsilon_{0}^{-4}\abs{\mI^{(\ell^{*})}(f)})$ \COMMENT{zero of $t\mapsto\Phi^{(\ell^{*})}(\omega_{k},t)$, see \eqref{eq:indicators-2}}
\ENDFOR 
\RETURN $D_{\rm rec} = \bigcap_{k}\{x\in B_{1} : \omega_{k}\cdot x\ge t_{\omega_{k}}^{(\ell^{*})}\}$
\end{algorithmic} 
\end{algorithm}

\begin{figure}[H]
    \centering
    \includegraphics[width=0.5\linewidth]{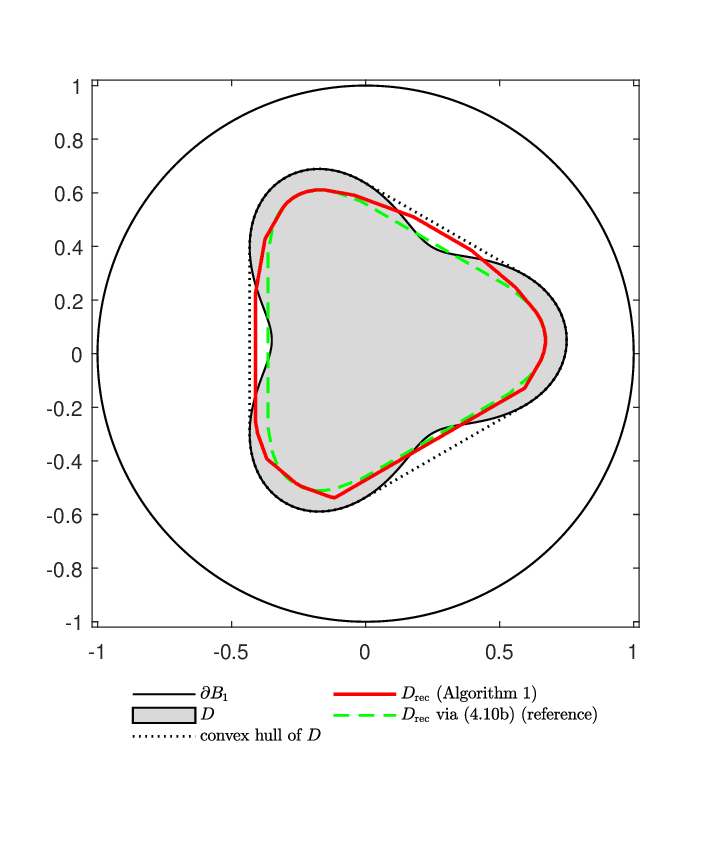} 
    \caption{Reconstruction of the convex hull of $D$ by \Cref{alg:implementation} (solid) and the reference via \eqref{eq:indicator-numerical-analysis-approx} (dashed), parameters as in \Cref{tab:param}}
    \label{fig:reconstruction}
\end{figure}

In line 7, we compute $W^{(\ell)}$ through the Dirichlet Green function of $B_{1}$, and the normal derivative in line 9 through the Poisson kernel $P(\hat{x},y)=(1-\abs{y}^{2})/(2\pi\abs{\hat{x}-y}^{2})$, using the Gauss-Legendre (radial) $\times$ trapezoidal (angular) quadrature with $40\times240$ nodes in $D$ and the trapezoidal rule with $1024$ nodes in $\partial B_{1}$. The normal derivative $\partial_{\nu}v_{f}$ is evaluated from its closed form. 
The level $t_{\omega}^{(\ell^{*})}$ in \Cref{alg:implementation} differs from its linearized counterpart 
\begin{equation*}
t_{\omega}^{(1)} := -1-\frac{h}{4}\log K^{(1)}(\omega) = -\frac{h}{4}\log\int_{D} e^{-\frac{4}{h}\omega\cdot x}\,\rmd x, 
\end{equation*}
by $-\frac{h}{4}\log(\abs{K^{(\ell^{*})}(\omega)}/K^{(1)}(\omega))$, which is negligible for small $h$ by \Cref{thm:approximation}\ref{itm:c}. 
Finally, we implement \Cref{alg:implementation} using the parameter choices listed in \Cref{tab:param}; the resulting reconstruction is shown in \Cref{fig:reconstruction}. 

\begin{table}[H]
\rowcolors{2}{white}{lightgray!50}
\renewcommand{\arraystretch}{1.5} 
\centering
\begin{tabular}{c|c}
parameter & value \\ 
\hline\hline 
$D$ & $\begin{aligned}
& \left\{(0.05,0.05) + sr(\varphi)(\cos\varphi,\sin\varphi) : 0\le s<1, 0\le\varphi<2\pi \right\} \\ 
& \quad r(\varphi)=0.55 + 0.15\cos(3\varphi)
\end{aligned}$ \\ 
$h$ & $0.05$ \\ 
$\varepsilon_{0}$ & $1/\sqrt{2}\approx0.707$ \\ 
$J$ & $J_{t}=t+1+h\log(1/\varepsilon_{0})\approx t+1.017$ \\ 
$\ell^{*}$ & $3$ \\ 
$N_{\omega}$ & $36$ 
\end{tabular}
\caption{Parameter settings}
\label{tab:param}
\end{table}

\Cref{fig:reconstruction} shows, as a dashed curve, the reconstruction obtained when the indicator is evaluated through the identity \eqref{eq:indicator-numerical-analysis-approx} instead of line~10 in \Cref{alg:implementation}. Since \eqref{eq:indicator-numerical-analysis-approx} uses $D$, this curve is not a reconstruction, it only indicates what \Cref{alg:implementation} would return in exact arithmetic. Both curves lie inside the convex hull of $D$.

\end{sloppypar}

\bibliographystyle{custom}
\bibliography{ref}
\end{document}